\documentclass[12pt]{amsart}
\usepackage{amssymb, mathrsfs,amsmath}
\usepackage{latexsym}
\usepackage{amsfonts}
\usepackage{shadow}

\newtheorem{Pa}{Paper}[section]
\newtheorem{Tm}[Pa]{{\bf Theorem}}
\newtheorem{La}[Pa]{{\bf Lemma}}
\newtheorem{Dn}[Pa]{{\bf Definition}}
\newtheorem{Cy}[Pa]{{\bf Corollary}}

\newtheorem{Rk}[Pa]{{\bf Remark}}
\newtheorem{Pn}[Pa]{{\bf Proposition}}

\newtheorem{Ex}[Pa]{{\bf Example}}
\newtheorem{Ap}[Pa]{{\bf Application}}

\newcommand{\w}{\omega}
\newcommand{\la}{\lambda}
\newcommand{\HH}{\mathcal H}

\date{}
\author[D. Alpay]{Daniel Alpay}
\author[P. Jorgensen]{Palle Jorgensen}
\address{(DA) Department of Mathematics \newline
Ben Gurion University of the Negev \newline P.O.B. 653, \newline
Be'er Sheva 84105, \newline ISRAEL} \email{dany@math.bgu.ac.il}
\address{(PJ)
Department of Mathematics\newline 14 MLH \newline The University
of Iowa, Iowa City,\newline IA 52242-1419 USA}
\email{jorgen@math.uiowa.edu}
\title[Spectral theory and generators for Gaussian processes]
{Spectral theory for Gaussian processes: Reproducing kernels,
random functions, boundaries, and $\mathbf L^2$-wavelet
generators with fractional scales} 
\dedicatory{To the Memory of William B. Arveson}
\begin{document}
\maketitle
\begin{abstract}
A recurrent theme in functional analysis is the interplay between
the theory of positive definite functions, and their reproducing
kernels, on the one hand, and Gaussian stochastic processes, on
the other. This central theme is motivated by a host of
applications, e.g., in mathematical physics, and in stochastic
differential equations, and their use in financial models. In
this paper, we show that, for three classes of cases in the
correspondence, it is possible to obtain explicit formulas which
are amenable to computations of the respective Gaussian
stochastic processes. For achieving this, we first develop two
functional analytic tools. They  are: $(i)$  an identification of
a universal sample space $\Omega$ where we may realize the
particular Gaussian processes in the correspondence; and (ii) a
procedure for discretizing computations in $\Omega$. The three
classes of processes we study are as follows: Processes
associated with: (a) arbitrarily given sigma finite regular
measures on a fixed Borel measure space; (b) with Hilbert spaces
of sigma-functions; and (c) with systems of self-similar measures
arising in the theory of iterated function systems. Even our
results in (a) go beyond what has been obtained previously, in
that earlier studies have focused on more narrow classes of
measures, typically Borel measures on $\mathbb R^n$. In our last
theorem (section 10), starting with a non-degenerate positive
definite function $K$  on some fixed set $T$, we show that there
is a choice of a universal sample space $\Omega$, which can be
realized as a "boundary" of $(T, K)$. Its boundary-theoretic
properties are analyzed, and we point out their relevance to the
study of electrical networks on countable infinite graphs.
\end{abstract}
\tableofcontents
\section{Introduction}
\setcounter{equation}{0}
We are considering three functional analytic questions arising at
the crossroads of pure and applied probability theory. In
different contexts of non-deterministic analysis, one needs
mathematical representations of the set of all possible outcomes,
called {\it the sample space} $\Omega$, of some experiment, for
example involving random trials.  This is easy enough in simple
discrete models, for example in experiment with tossing coins.
The sample space of each trial is the set \{head, tail\}, and
more subtle models then involve Cartesian products. However in
infinite models, and in most continuous models, a complete
description of a sample space of outcomes and its subsets,
events, presents subtle problems. In Brownian motion models, for
example, $\Omega$ may conveniently be taken to be a suitable
space of continuous functions, sample paths. Now, to approach
computations, one is faced with the use of simulations of
suitable subsets in $\Omega$; e.g., Monte-Carlo simulations. For
such approaches, because of noise, of uncertainties, or limited
information, it is often helpful to pick different mathematical
realizations of the set $\Omega$: For example, a version of
$\Omega$ consisting of sample paths defined only on suitable
subsets, as opposed to defined point-wise. This is often good
enough as one is interested in particular functions on $\Omega$.
Whichever choice is made, $\Omega$ will naturally come equipped
with a sigma-algebra, say $\mathscr F$, of subsets, and a
probability measure $P$ defined on $\mathscr F$. The $\mathscr F$
-measurable functions are random variables, and systems of random
variables are stochastic processes. The process is Gaussian if we
can choose the probability measure $P$ such that the random
variables making up the process are Gaussian, and in $\mathbf
L^2(P)$.\\

With the use of the corresponding Gaussian densities, and
covariance functions, one then computes quantities from the
random variables; and the question of choice of $\Omega$ can then
often be avoided. Nonetheless, for applications to stochastic
integration, one is forced to be more precise with the choice of
$\Omega$, and a number of functional analytic tools are available
for the purpose. In the approach to this problem based on
Gelfand-triples (see Section \ref{sec5}), one may realize $\Omega$
as a space of Schwartz-tempered distributions. However with this
realization of $\Omega$, it is more difficult to make a direct
connection to the initial model, and to set up suitable
Monte-Carlo simulations. As a result, there is a need for
discretizations. Several such discretizations will be presented
here, and comparisons will be made.\\

Our approach, in this general context, relies on our use of Gaussian
Hilbert spaces, and of associated sequences of independent,
identically distributed (i. i. d.) standard Gaussian $N(0,
1)$-random variables. But this then further introduces a host of
choices, and of these we identify one which is universal in a
sense made precise in Sections \ref{sec5}-\ref{sec7}.\\

In this paper, we will focus on Gaussian stochastic processes,
but we also offer applications of our results to certain random
functions (Section \ref{sec9}) which involve non-Gaussian
distributions. Similarly, a host of
simulation approaches involve non-Gaussian choices.\\

The purpose of the paper is three-fold. First we study $(i)$ a
universal choice of sample space for a family of $\mathbf L^2$
Gaussian noise processes. While these processes have appeared in
one form or the other in prior literature, the choice of sample
spaces has not been studied in a way that facilitates comparisons.
We index these Gaussian noise processes by the set of regular
measures in some fixed measure space $(M, \mathscr B)$, with
$\mathscr B$ some given Borel sigma-algebra of subsets in M.
Secondly we make precise $(ii)$ equivalence in this category of
Gaussian noise processes, and we prove a uniqueness theorem, where
uniqueness is specified by a specific measure isomorphism of the
respective sample spaces. In our third result $(iii)$, given a
fixed measure space $(M, \mathscr B)$, we identify a Hilbert space
$\HH$ , with the property that the Gaussian noise process indexed
by $\HH$ universal envelope of all the Gaussian noise processes
from $(ii)$. As applications we compute Gaussian noise processes
associated to Cantor measures, and more generally to iterated
function systems (IFS) measures, and to a family of
reproducing kernel Hilbert spaces (RKHS).\\

For readers not familiar with Gaussian processes, for the present
purpose, the following are helpful:
\cite{aal2,MR2793121,alp,MR0265548,Hida_BM,MR2083706,MR1203453,MR0102759,MR0264761};
for infinite products and applications, see \cite{Ka48},
\cite{MR1557013}, and  \cite{aj-cras}, \cite{ajlm2}. The
universal Hilbert space from (iii) is used in a different context
\cite{MR1978577,Hida_BM,MR2097020,MR0214150,Nelson_flows}. For a
small sample of recent applications, we cite
\cite{MR2133804,MR2444857,MR1960424} For reproducing kernel
Hilbert spaces, see, for example, \cite{MR2002b:47144,saitoh}. In
the way of presentation, it will be convenient to begin with a
quick review of infinite products, this much inspired by the
pioneering paper \cite{Ka48} by Kakutani.

\section{Preliminaries}
\setcounter{equation}{0} \label{sec2} Below we present a framework
of Gaussian Hilbert spaces. These in turn play a crucial role in
the study of positive semi-definite kernels, and their associated
reproducing kernel Hilbert spaces, see Sections
\ref{sec9}-\ref{sec10}. In its most general form, the theory of
Gaussian Hilbert spaces $\HH $ is somewhat abstract, and it is
therefore of interest, for particular cases of $\HH$, to study
natural decompositions into cyclic components in $\HH$ which
arise in applications, and admit computation. Hence we begin with
those processes whose covariance function may be determined by a
fixed measure. Even this simpler case generalizes a host of
Gaussian processes studied earlier with the use of Gelfand
triples built over the standard Hilbert space $\mathbf L^2(
\mathbb R^d,dx)$, with $dx$ denoting the Lebesgue measure, with
the use of Laurent Schwartz theory of tempered distributions. Our
present framework is not confined to the Euclidean case. Indeed,
starting with any measure space $M$ and a Borel sigma-algebra
$\mathscr B$, we then show in Section \ref{sec4} that the General
Gaussian Hilbert space (Definition \ref{def2.2}) decomposes as an
orthogonal sum where the corresponding cyclic subspaces are those
generated by a family of sigma-finite measures on $M$.  Indeed,
in applications to measurement, in physics, and in statistics, it
is often not possible to pin down a variable as a function of
points in the underlying
space $M$. As a result, it has proved useful to study processes indexed by sigma-algebras of subsets of $M$.\\

In our consideration of random variables, of Hilbert spaces, and
of Gaussian stochastic processes, it will be convenient for us to
restrict to the case of {\it real-valued} functions and real
Hilbert spaces. It will be helpful to first state the respective
results in the real case, and then, at the end, when needed,
remove the restriction. One instance when complex Hilbert spaces
are needed is the introduction of Fourier bases, i.e., orthogonal
bases consisting of functions $e_\lambda$ where
$\lambda\in\mathbb R$ and $e_\lambda(x)=e^{i\lambda x}$ or
$e^{2\pi i \lambda x}$. However, our setting will be general
measure spaces $(M,\mathscr B,\mu)$, where $\mathscr B$ is a
sigma-algebra of measurable subsets of some set $M$, and $\mu$ is
a positive measure on $M$. The restricting assumption is {\sl
sigma-finiteness}, i.e., there are subsets $B_1,B_2,\ldots$ of
$\mathscr B$  such that
\begin{equation}
\label{eq:1}
M=\bigcup_{j=1}^\infty B_j,\quad{\rm and}\quad
\mu(B_j)<\infty,\,\,\,\forall j\in\mathbb N.
\end{equation}

\begin{Dn}
\label{def2.1}
A Gaussian (noise) stochastic process indexed by
$(M,\mathscr B,\mu)$ consists of a probability space
$(\Omega,\mathscr F,\mathbb P)$: $\Omega$ is a set (sample
space), $\mathscr F$ is a sigma-algebra of subsets (events) of
$\Omega$, and $\mathbb P$ is a probability measure defined on
$\mathscr F$. We assume that, for all $A\in \mathscr B$ such that
$\mu(A)<\infty$, there is a Gaussian random variable
\begin{equation}
W_A=W_A^{(\mu)}\,\,:\,\,\Omega\longrightarrow\mathbb R
\end{equation}
with zero mean and variance $\mu(A)$ (that is, $W_A\sim
N(0,\mu(A))$, the Gaussian with zero mean and variance $\mu(A)$),
i.e. for all $a,b\in\mathbb R\cup\left\{\pm \infty\right\}$ with
$a<b$,
\[
\left\{\w\in\Omega\,\,\,|\,\,\, a<W_A(\w)\le b\right\}\in \mathscr F,
\]
and
\[
\begin{split}
\mathbb P\left(\left\{a<W_A(\w)\le b\right\}\right)&=\int_a^b\frac{1}{\sqrt{2\pi \mu(A)}}e^{-\frac{x^2}{2\mu(A)}}dx\\
&=\int_{\frac{a}{\sqrt{\mu(A)}}}^{\frac{b}{\sqrt{\mu(A)}}}\frac{1}{\sqrt{2\pi \mu(A)}}e^{-\frac{x^2}{2}}dx\\
&=\gamma_1((\frac{a}{\sqrt{\mu(A)}},\frac{b}{\sqrt{\mu(A)}}]),
\end{split}
\]
where $\gamma_1$ is the standard Gaussian on $\mathbb R$.
\end{Dn}

\begin{Dn}
\label{def2.2} A Gaussian process indexed by a (fixed) Hilbert
space $\mathcal H$ consists of a probability space
$(\Omega,\mathscr F,\mathbb P)$ such that, for all $F\in\mathcal
H$, there is a Gaussian random variable $W_F$ with law
$N(0,\|F\|^2_{\mathcal H})$ such that
\begin{equation}
\mathbb E(W_{F_1}W_{F_2})=\langle F_1,F_2\rangle_{\mathcal H},
\quad \forall F_1,F_2\in\mathcal H.
\end{equation}
It is further assumed that for all
$\left\{A_j\right\}_{j=1}^n\subset\mathscr B$ such that
$0<\mu(A_j)<\infty$, $i=1,2,\ldots, n$, the joint distribution of
the family $\left\{W_{A_j}\right\}_{j=1}^n$ is Gaussian with zero
mean and covariance matrix $(\mu(A_i\cap A_j))_{i,j=1}^n$. We
will assume throughout that ${\rm dim}~\mathbf L^2(\mu)=\infty$.
The finite dimensional case is dealt with separately.
\end{Dn}

\begin{Rk}
\label{rk2.3}
With the specifications in Definitions \ref{def2.1} and \ref{def2.2}, it is known that, in each case, such Gaussian
processes exist. In the case of Definition \ref{def2.2}, when $\left\{F_i\right\}_{i=1}^n$ is a system
in $\mathcal H$, then the random variables
$\left\{W_{F_i}\right\}_{i=1}^n$ have a joint Gaussian
distribution corresponding to the covariance matrix
$\left(\langle F_i,F_j\rangle_{\mathcal H}\right)_{i,j=1}^n$.
\end{Rk}

Starting with a measure space $(M, \mathscr B)$ , we will show in
Section \ref{sec4}, that there is a universal Hilbert space
$\mathcal H$ which "contains" all the stochastic processes
derived from sigma-finite measures  $\mu$ on $(M, \mathscr B)$. In
detail, given an arbitrary  $\mu$ , we get a Gaussian process
$W^{(\mu)}$ with $\mu$ as its covariance measure; see Definition
\ref{def2.1}. Now, the universal Hilbert space $\mathcal H$ over
$(M, \mathscr B)$ will satisfy the conditions in Definition
\ref{def2.2}; and it will be a Hilbert space of sigma-functions
(Definition \ref{def3.1}). Before getting to this, we must prepare
the ground with some technical tools. This is the purpose of the
next section on infinite products, and discrete Gelfand-triples.

\section{The probability space $(\Omega_{\mathbf s},\mathscr F_{\mathbf s},Q)$}
\setcounter{equation}{0} \label{sec5}
The purpose of this section is to show that there is a single
infinite-product measure space such that for every measure space
$M$ and fixed Borel sigma-algebra $\mathscr B$, everyone of the
Gaussian processes $W^{(\mu)}$, where $\mu$ is sigma-finite
measure on $M$, may be represented in $\mathbf L^2$ of this
infinite-product measure space. Since the construction must apply
to every sigma-finite measure $\mu$, we must adjust the
construction so that it can be adapted to orthonormal bases
(ONBs) in each of the corresponding $\mathbf L^2(\mu)$ Hilbert
spaces. To do this, we will be introducing a suitable Gelfand
triple (see \eqref{eq28}-\eqref{eq30}), realized in sequence
spaces, as opposed to the more traditional setting based instead
on $\mathbf L^2( R^d, dx)$ and Schwartz' tempered
distributions. There is a number of advantages of this approach,
for example we are not singling out any particular $\mathbf
L^2(\mu)$, and also not a
particular choice of ONB.\\

An initial choice for $\Omega_S$ is $\Omega_S=\times_{\mathbb
N}\mathbb R$, that is the space of all functions from $\mathbb N$
into $\mathbb R$, or equivalently, of all real sequences
$(c_1,c_2,\ldots)$ indexed by $\mathbb N$. Let $\mathbf s$ be the
space of sequences $c=(c_n)_{n\in\mathbb N}\in\times_{\mathbb
N}\mathbb R$ with the following property: For every $p\in\mathbb
N$ there exists $K_p<\infty$ such that
\begin{equation}
\label{27} |c_j|\le K_p j^{-p},\quad \forall j\in\mathbb N,
\end{equation}
and denote by $\mathbf s^\prime$ the dual space of all sequences
$\xi=(\xi_j)_{j\in\mathbb N}$ of polynomial growth, that is, such
that there exists $q\in\mathbb N$ and $K_q>0$ such that
\begin{equation}
\label{eq28} |\xi_j|\le K_q j^q,\quad \forall j\in\mathbb N.
\end{equation}
Then (see \cite{MR2083706,MR1408433,MR751959})
\[
\mathbf s\subset \ell^2\subset \mathbf s^\prime
\]
is a Gelfand triple, i.e., with the semi-norms defined from
\eqref{27}, $\mathbf s$ becomes a Fr\'echet space, and the
embedding from $\mathbf s$ into $\ell^2$ is nuclear (and $\mathbf
s^\prime$ denotes the dual of $\mathbf s$).\\

Let $\mathscr F_{\mathbf s}$ denote the sigma-algebra of subsets
in $\mathbf s^\prime$ generated by the cylinder sets as follows:
For $c_1,c_2,\ldots, c_n\in\mathbf s$ and an open set
$O\subset\mathbb R^n$, define the cylinder ${\rm
Cyl}~(c_1,\ldots,c_n,O)$ by
\begin{equation}
\label{eq30} {\rm Cyl}~(c_1,\ldots, c_n,O)=\left\{\xi\in \mathbf
s^\prime\,\, |\,\, \left(\langle \xi, c_1\rangle,\ldots, \langle
\xi, c_n\rangle\right)\in O\right\}.
\end{equation}
As the data in \eqref{eq30} varies, we get the cylinder sets in
$\mathbf s^\prime$ and the corresponding
sigma-algebra $\mathscr F_{\mathbf s}$.\\

Further note that the sets in \eqref{eq30} generate a system of
neighborhoods for the ${\rm weak}^*$-topology on $\mathbf
s^\prime$. Moreover, if $\mathbf s$  is assigned its Fr\'echet
topology from the semi-norms in \eqref{27}, then $\mathbf
s^\prime$ (with its ${\rm weak}^*$-topology) is the dual of
$\mathbf s$.

\begin{La}
\label{le5.1}
From Gelfand's theory we therefore get the
existence of a unique probability measure $Q$ on $(\mathbf
s^\prime, \mathscr F_{\mathbf s})$ with the property that
\begin{equation}
\label{eq31} \int_{\mathbf s^\prime} e^{i\langle
\xi,c\rangle}dQ(\xi)=e^{-\frac{1}{2}\|c\|^2_2},
\end{equation}
where
\begin{equation}
\label{eq321} \langle \xi, c\rangle=\sum_{j=}^\infty c_j\xi_j,
\end{equation}
and $\|c\|^2_2=\sum_{j=1}^\infty c_j^2$. Moreover, $\langle
\xi,c\rangle$ in \eqref{eq31} and \eqref{eq321} extends from
$\mathbf s\times \mathbf s^\prime$ to $\ell^2\times \mathbf
s^\prime$, representing every $c\in\ell^2$ as a Gaussian variable
on $(\mathbf s^\prime,\mathscr F_{\mathbf s},Q)$, with
\begin{equation}
\label{eq32} \mathbf E_Q\left(\langle \cdot, c\rangle\right)=0,
\end{equation}
and
\begin{equation}
\label{eq33} \mathbf E_Q\left(\langle \cdot,
c\rangle^2\right)=\|c\|^2_2.
\end{equation}
Furthermore, the set of coordinate functions on $\Omega_{\mathbb
S}:\mathbf s^\prime$,
\[
\pi_j(\xi)=\xi_j,\quad j\in\mathbb N,
\]
turns into an independent, identically distributed (i.i.d.)
system of $N(0,1)$ standard Gaussian variables, and we get:
\begin{equation}
\label{eq34} \mathbf E\left(\pi_{j_1}\pi_{j_2}\cdots
\pi_{j_k}e^{i\langle \cdot,
c\rangle}\right)=(-1)^{k/2}c_{j_1}c_{j_2}\cdots c_{j_k}
e^{-\frac{1}{2}\|c\|_2^2}.
\end{equation}
\end{La}

{\bf Proof:} We begin with the assertion
\begin{equation}
\label{eq508}
\mathbb E_Q\left(\pi_j\pi_k\right)=\delta_{j,k},\quad\forall
j,k\in\mathbb N.
\end{equation}
Take first $j=k$; then,
\[
\mathbb E_Q\left(\pi_j^2\right)=\frac{1}{\sqrt{2\pi}}\int_{\mathbb
R}c_j^2e^{-\frac{c_j^2}{2}}dc_j=1,
\]
and if $j\not =k$ we get
\[
\begin{split}
\mathbb E_Q\left(\pi_j\pi_k\right)&=\frac{1}{2\pi}\iint_{\mathbb
R^2}c_jc_ke^{-\frac{c_j^2+c_k^2}{2}}dc_jdc_k\\
&=\left(\frac{1}{\sqrt{2\pi}}\int_{\mathbb
R}xe^{-\frac{x^2}{2}}dx\right)^2\\
&=0,
\end{split}
\]
which proves \eqref{eq508}.\\

We now prove the assertion \eqref{eq31}. With \eqref{eq321} we get
\[
\langle\xi,c\rangle=\sum_{j=1}^\infty c_j\pi_j(\xi),\quad\forall
\xi\in\mathbf s^\prime,\,\,{\rm and}\,\,\forall c\in\mathbf s.
\]
We prove \eqref{eq31} for $c\in\mathbf s$, and then extend it to
all of $\ell^2$. We have
\[
\begin{split}
\int_{\mathbf s^\prime} e^{i\langle \xi,c\rangle}dQ(\xi)
&=\mathbb E_Q\left(e^{i\sum_{k=1}^\infty
c_k\pi_k(\cdot)}\right)\\
&=\prod_{k=1}^\infty \mathbb E_Q\left(e^{ic_k\pi_k(\cdot)}\right)\\
&=\prod_{k=1}^\infty e^{-\frac{c_k^2}{2}}\\
&=e^{-\frac{1}{2}\|c\|^2_2},
\end{split}
\]
which is the desired conclusion.\\

The proof of \eqref{eq34} follows from an application of
\eqref{eq31} to
\begin{equation}
\label{eq1000}
c+t_1e_{j_1}+\cdots +t_ke_{j_k},\quad t_1,\ldots, t_k\in\mathbb R,
\end{equation}
where
\[
(e_j)_\ell:=\delta_{j,\ell},\quad\forall j,\ell\in\mathbb N.
\]
is the standard ONB in $\ell^2$, i.e. $\widetilde{e_j}=\pi_j$.
Now \eqref{eq34}  follows if \eqref{eq1000} is substituted into
\eqref{eq31}, and the partial derivatives
$\frac{\partial}{\partial t_1}\cdots\frac{\partial}{\partial t_k}$
are computed on both sides, and then evaluated at $t_1=\cdots
=t_k=0$.
\mbox{}\qed\mbox{}\\

\begin{La}
\label{la1001}
Let $Q=\times_{\mathbb N}\gamma_1$ be the product measure.  Then
\begin{equation}
\label{eq1003} Q(\mathbf s^\prime)=1\quad and\quad Q(\ell^2)=0.
\end{equation}
\end{La}

{\bf Proof:} The first claim follows from Minlos's theorem
applied to the following positive definite function on $\mathbf s$
\begin{equation}
\label{eq1002}
c\in\mathbf s\quad\mapsto\quad e^{-\frac{\|c\|^2}{2}}.
\end{equation}
Indeed, the function \eqref{eq1002} is clearly continuous with
respect to the semi-norms in $\mathbf s$; see \eqref{27}.\\

We need to prove the second claim in \eqref{eq1003}, i.e. the
assertion that $\ell^2\subset \mathbf s^\prime$ has $Q$-measure
zero. Assume the contrary, i.e.  assume $Q(\ell^2)>0$. Since
\[
\lim_{j\rightarrow\infty}\pi_j=0
\]
point-wise on $\ell^2$, we have
\[
\lim_{j\rightarrow\infty}\int_{\ell^2}e^{i\pi_j(\w)}dQ(\w)=Q(\ell^2)
\]
as an application of Lebesgue's dominated convergence theorem. On
the other hand,
\[
\mathbb E_Q(e^{i\pi_k(\cdot)})=e^{-\frac{1}{2}},\quad\forall
k\in\mathbb N,
\]
and so another application of Lebesgue's dominated convergence
theorem leads to
\[
e^{-\frac{1}{2}}=Q(\ell^2)+Q(\mathbf s^\prime\setminus\ell^2),
\]
which is a contradiction since the sum should be equal to $1$.
\mbox{}\qed\mbox{}\\

\begin{Tm}
\label{tm5.1}
Let $(M,\mathscr B,\mu)$ be a sigma-finite measure
space as specified in Section \ref{sec3}, and let
$\left\{\varphi_j\right\}_{j\in\mathbb N}$ be a choice of
orthonormal basis (ONB) in $\mathbf L^2(\mu)$. Then, the Gaussian
process $W^{(\mu)}$ may be realized in $\mathbf
L^2(\Omega_{\mathbb S},\mathscr F_{\mathbf s},Q)$ as follows: For
$A\in\mathscr B$ such that $0<\mu(A)<\infty$, set
\begin{equation}
\label{eq35}
W_A^{(\mu)}(\xi)=\sum_{j=1}^\infty\left(\int_A\varphi_j(x)d\mu(x)\right)
\pi_j(\xi),\quad \xi\in{\mathbf s}^\prime.
\end{equation}
Then, $W^{(\mu)}$, defined by \eqref{eq35}, is a copy of the
Gaussian process from Definition \ref{def2.1}.
\end{Tm}

{\bf Proof:} In view of \eqref{eq34}, we need only to prove that
$W^{(\mu)}_A$ in \eqref{eq35} is a $N(0,\mu(A))$ Gaussian
variable, and that
\begin{equation}
\label{eq37} \mathbb
E_Q\left(W_A^{(\mu)}W_B^{(\mu)}\right)=\mu(A\cap B), \quad
\forall A,B\in\mathscr B.
\end{equation}
But the first assertion follows from
\[
\sum_{j=1}^\infty\left(\int_A\varphi_j(x)d\mu(x)\right)^2=
\|\chi_A\|^2_{\mathbf L^2(\mu)}=\mu(A),
\]
and we prove \eqref{eq37} as follows:
\[
\begin{split}
\mathbb
E_Q\left(\left(\sum_{j=1}^\infty\left(\int_A\varphi_j(x)d\mu(x)\right)
\pi_j\right)
\left(\sum_{j=1}^\infty\left(\int_B\varphi_k(x)d\mu(x)\right)\pi_k\right)\right)&=\\
&\hspace{-7cm}=\sum_{j,k=1}^\infty\left(\int_A\varphi_j(x)d\mu(x)
\right)\left(\int_B\varphi_k(x)d\mu(x)\right)\delta_{j,k}\\
&\hspace{-7cm}=\langle \chi_A,\chi_B\rangle_{\mathbf L^2(\mu)}\\
&\hspace{-7cm}=\mu(A\cap B),
\end{split}
\]
which is the desired conclusion.
\mbox{}\qed\mbox{}\\

In the next section, we generalize the expansion formula
\eqref{eq35} above.

\begin{Cy}
\label{cy5.4}
Let $(M,\mathscr B)$ be as in Theorem \ref{tm5.1}, let $\mu$ and
$\lambda$ be two sigma-finite measures defined on it, such that
$\mu<<\lambda$, and let $f\in\mathbf L^2(M,\mathscr M,\mu)$.
Then, in the representation \eqref{eq35}, referring to $\mathbf
L^2(\mathbf s^\prime, Q)$, we have
\begin{equation}
\label{eq5.18}
W^{(\lambda)}\left(f\sqrt{\frac{d\mu}{d\lambda}}\right)=
W^{(\mu)}(f),\quad Q\,\, a.e.
\end{equation}
\end{Cy}

{\bf Proof:} Picking an ONB $\left\{\varphi_j\right\}_{j\in\mathbb
N}$ in $\mathbf L^2(M,\mathscr M,\mu)$, we note that then
$\left\{\varphi_j\sqrt{\frac{d\mu}{d\lambda}}\right\}_{j\in\mathbb
N}$ is an ONB in $\mathbf L^2(M,\mathscr M,\lambda)$. Now use
\eqref{eq35} for the pair of ONBs. We get
\[
\begin{split}
W^{(\mu)}(f)&=\sum_{j=1}^\infty\langle\varphi_j,f\rangle_{\mathbf
L^2(\mu)}\pi_j\\
&=\sum_{j=1}^\infty
\left(\int_M\varphi_j(x)f(x)d\mu(x)\right)\pi_j\\
&= \sum_{j=1}^\infty
\left(\int_M\varphi_j(x)f(x){\frac{d\mu}{d\lambda}(x)}d\lambda(x)\right)\pi_j\\
&= \sum_{j=1}^\infty \left(\int_M \sqrt{\frac{d\mu}{d\lambda}(x)}
\varphi_j(x)\sqrt{\frac{d\mu}{d\lambda}(x)}f(x)d\lambda(x)\right)\pi_j\\
&=W^{(\lambda)}\left(f\sqrt{\frac{d\mu}{d\lambda}}\right).
\end{split}
\]

\mbox{}\qed\mbox{}\\

We need another construction of a universal space as well, using a
construction of Kakutani \cite{Ka48}. More precisely, consider the
space $\times_{\mathbb N}\mathbb R$, and denote by $\xi$ a running
element in this cartesian product. Define for
$F(\xi)=f_n(\xi_1,\ldots ,\xi_n)$, where $f_n$ is a measurable and
summable function of $n$ real variables
\[
\mathcal L(F)=\int\cdots\int_{\mathbb R^n}f_n(\xi_1,\ldots,
\xi_n)\gamma_n(\xi_1,\ldots ,\xi_n)d\xi_1\cdots d\xi_n,
\]
where $\gamma_n$ is the product of the densities of $n$ i.i.d.
$N(0,1)$ variables. By Kolmogorov's theorem \cite{MR0622034},
there exists a unique probability $Q_K$ on $\prod_{\mathbb
N}\mathbb R$ such that
\[
\int_{\prod_{\mathbb N}\mathbb
R}F(\xi)dQ_K(\xi)=\int\cdots\int_{\mathbb R^n}f_n(\xi_1,\ldots,
\xi_n)\gamma_n(\xi_1,\ldots ,\xi_n)d\xi_1\cdots d\xi_n.
\]
In fact,
\begin{equation}
\label{eq18}
Q_K=\prod_{\mathbb N}\gamma_1
\end{equation}
on the countably infinite Cartesian product $\prod_{\mathbb N}
\mathbb R$, where $\gamma_1$ is the standard $N(0,1)$ Gaussian on
$\mathbb R$ (with density
$\frac{1}{\sqrt{2\pi}}e^{-\frac{x^2}{2}}$).\\

The measure $Q_K$ and $Q$ have the same characteristic function,
and $Q_K(\mathbf s^\prime)=1$. So we will in the sequel use both
the spaces $(\mathbf s^\prime,\mathscr F_{\mathbf s}, Q)$ and
$(\prod_{\mathbb N}\mathbb R, Q_K)$.

%Unless specified otherwise, we will now assume that all the
%$\mathbf L^2$-Gaussian processes are realized in the space
%$\mathbf L^2(\Omega_s,\mathscr F_s,Q)$. The associated
%expectation will be denoted by $\mathbb E_Q$.

\section{A Hilbert space of sigma-functions}
\setcounter{equation}{0}
\label{sec3}
In spectral theory, in
representation theory (see e.g. \cite{MR1978577,Nelson_flows}),
and in the study of infinite products \cite{Hida_BM}, and of
iterated function systems (IFS) (see e.g. \cite{MR625600}) one is
faced with the problem of identifying direct integral
decompositions. Naturally, a given practical problem may not by
itself entail a Hilbert space, and, as a result, one must be
built by use of the inherent geometric features of the problem.
In these applications it has proved useful to build the Hilbert
space from a set of equivalence classes. The starting point will
be pairs  $(f, \mu)$ where $\mu$ is a measure, and $f$ is a
function, assumed in $\mathbf L^2(\mu)$. It turns out (see
\cite{Nelson_flows}) that the set of such equivalence classes
acquire the structure of a Hilbert space, called a sigma-Hilbert
space. Further we show through applications (Sections \ref{sec7}
and \ref{sec8})  that these sigma-Hilbert spaces form a versatile
tool in the study of Gaussian processes. These Gaussian processes
are indexed by a choice of a suitable sigma-algebras of subsets
of $M$.

\begin{Dn}
\label{def3.1}
Let $(M,\mathscr B)$ be a fixed measure space, and let $\mathscr
M(M,\mathscr B)$ denote the set of all sigma-finite positive
measures on $(M,\mathscr B)$. For pairs $(f_i,\mu_i)$, $i=1,2$,
where
\begin{equation}
\label{eq6}
\mu_i\in\mathscr M(M,\mathscr B),\quad f_i\in \mathbf L^2(\mu_i),
\end{equation}
we introduce the equivalence relation $\sim$ as follows:
$(f_1,\mu_1)\sim (f_2,\mu_2)$ if and only if there exists
$\lambda\in\mathscr M(M,\mathscr B)$ such that $\mu_i<<\lambda$
and
\begin{equation}
\label{eq7}
f_1\sqrt{\frac{d\mu_1}{d\lambda}}=f_2\sqrt{\frac{d\mu_2}{d\lambda}},\quad
\lambda\,\, a.e.
\end{equation}
Here, $\frac{d\mu_i}{d\lambda}$ denote the respective
Radon-Nikodym derivatives. For the measure $\lambda$, we may take
$\lambda=\mu_1+\mu_2$.
\end{Dn}

It is known (see \cite{Nelson_flows}), that \eqref{eq7}
indeed defines an equivalence relation in the set of all pairs as
specified in \eqref{eq6}. If $\mu\in\mathscr M(M,\mathscr B)$ and
$f\in\mathbf L^2(\mu)$, we denote the equivalence class of
$(f,\mu)$ by $f\sqrt{d\mu}$. Moreover (see \cite{Nelson_flows}), set
\begin{equation}
\label{eq9}
\langle f_1\sqrt{d\mu_1},f_2\sqrt{d\mu_2}\rangle=\int_Mf_1(x)f_2(x)
\sqrt{\frac{d\mu_1}{d\lambda}(x)
\frac{d\mu_2}{d\lambda}(x)}d\lambda(x),
\end{equation}
where $\lambda$ is chosen such that $\mu_i<<\lambda$ for $i=1,2$ (for example,
one can take
$\lambda=\mu_1+\mu_2$) and set
\begin{equation}
\label{eq10}
f_1\sqrt{d\mu_1}+f_2\sqrt{d\mu_2}=
\mbox{{\rm equivalence class of}}~(f_1\sqrt{\frac{d\mu_1}{d\lambda}}+
f_2\sqrt{\frac{d\mu_2}
{d\lambda}},\lambda).
\end{equation}
The operations defined in \eqref{eq9} and \eqref{eq10} are known
to respect the equivalence relation \eqref{eq7}. The set of all
corresponding equivalence classes becomes a Hilbert space, which
we shall denote $\mathcal H=\mathcal H(M,\mathscr B)$. A separate
argument is needed in proving completeness, see
\cite{Nelson_flows}: If $(f_n\sqrt{d\mu_n})_{n\in\mathbb N}$ is a
Cauchy sequence in $\mathcal H(M,\mathscr B)$, there is a pair
$(f,\mu)$ with associated equivalence class $f\sqrt{d\mu}$ such
that
\[
\lim_{n\rightarrow\infty}\|f\sqrt{d\mu}-f_n\sqrt{d\mu_n}
\|_{\mathcal H(M,\mathscr B)}=0.
\]
\begin{Pn}
\label{rk3.2}
The Gaussian processes from Definition
\ref{def2.1}, with $(M,\mathscr B,\mu)$ given, are special cases
of the one from Definition \ref{def2.2} if we take $\mathcal
H=\mathbf L^2(\mu)$.
\end{Pn}
{\bf Proof:} To see this, fix $(M,\mathscr B,\mu)$, and let
$W^{(\mu)}$ be the associated Gaussian process. Then the map
\begin{equation}
\label{eq11} A\in\mathscr B,\,\,
\mu(A)<\infty\quad\Longrightarrow\quad W_A^{(\mu)}\in\mathbf
L^2(\Omega,\mathbb P)
\end{equation}
extends to all of $\mathbf L^2(\mu)$. The extended map, denoted by
\begin{equation}
\label{eq12}
W^{(\mu)}(f)=\int_Mf(x)dW_x^{(\mu)},
\end{equation}
and with range in $\mathbf L^2(\Omega,\mathbb P)$, is the Ito
integral \cite{Hida_BM}. When $f\in\mathbf L^2(\mu)$ is a simple
function, that is a finite sum of the form
\begin{equation}
\label{eq13}
f=\sum_{i=1}^N a_i\chi_{A_i},
\end{equation}
where the $a_i$ are real numbers, and the $A_i$ belong to
$\mathscr B$ are such that $A_i\cap A_j=\emptyset$ for $i\not
=j$, then set
\begin{equation}
W^{(\mu)}(f)=\sum_{i=1}^N a_iW_{A_i}^{(\mu)}.
\end{equation}
Using
\begin{equation}
\label{eq15}
\mathbb E(W_A^{(\mu)}W_B^{(\mu)})=\mu(A\cap B),\quad\forall A,B\in\mathscr B,
\end{equation}
we get
\begin{equation}
\label{eq16}
\mathbb E\left(|W^{(\mu)}(f)|^2\right)=\sum_{i=1}^Na_i^2\mu(A_i)=\int_Mf(x)^2d\mu(x).
\end{equation}
(In the complex case, we use $\int_M |f(x)|^2d\mu(x)$ on the
right hand-side of \eqref{eq16}). Since every function
$f\in\mathbf L^2(\mu)$ is the limit (in the norm of $\mathbf
L^2(\mu)$) of a sequence of simple functions, we conclude that
the isometry \eqref{eq16} extends to all of $\mathbf L^2(\mu)$.
Furthermore, by polarization,
\begin{equation}
\label{eq17}
\mathbb E\left(W^{(\mu)}(f)W^{(\mu)}(g)\right)=\langle f,g\rangle_{\mathbf L^2(\mu)},\quad\forall f,g\in\mathbf L^2(\mu).
\end{equation}
\mbox{}\qed\mbox{}\\

%\begin{Pn}
%\label {la3.3}
%There is a single probability space
%$(\Omega_s,\mathscr F_s,Q)$ such that all the Gaussian processes described above
%may be realized in the space $\mathbf L^2(\Omega_s,\mathscr F_s,Q)$.
%\end{Pn}

%{\bf Proof:} The details will be postponed, but

\begin{La}
\label{la3.4} Consider $(M,\mathscr B)$ as in Definition
\ref{def3.1}, and let $(f_i,\mu_i)$, $i=1,2$ be a pair, see
\eqref{eq6}. Then,
\begin{equation}
\label{eq19}
(f_1,\mu_1)\,\,\sim \,\, (f_2,\mu_2)\quad \iff\quad W^{(\mu_1)}(f_1)=W^{(\mu_2)}(f_2),\quad Q\,\, a.e.
\end{equation}
\end{La}

{\bf Proof:} We first assume that $(f_1,\mu_1)\,\,\sim \,\, (f_2,\mu_2)$. There exists $\lambda\in\mathscr M(M,
\mathscr B)$ such that both $\mu_1$ and $\mu_2$ are absolutely continuous with respect to $\lambda$ and such that
\eqref{eq7} is in force. Then,
\begin{eqnarray}
\label{eq3456}
W^{(\mu_1)}(f_1)&=&W^{(\lambda)}\left(f_1\sqrt{\frac{d\mu_1}{d\lambda}}\right)\\
\nonumber
&=&W^{(\lambda)}\left(f_2\sqrt{\frac{d\mu_2}{d\lambda}}\right)\\
\nonumber &=&W^{(\mu_2)}(f_2),
\end{eqnarray}
which is the desired identity on the right hand-side of
\eqref{eq19}. For the justification of \eqref{eq3456}, see
Section \ref{sec5}, especially
Corollary \ref{cy5.4}.\\

Conversely, assume that $W^{(\mu_1)}(f_1)=W^{(\mu_2)}(f_2)$ (almost everywhere with respect to $Q$) for some pairs
$(f_1,\mu_1)$ and $(f_2,\mu_2)$. By the argument above applied to $\lambda=\mu_1+\mu_2$, we get
\begin{equation}
\label{eq21}
W^{(\lambda)}\left(f_1\sqrt{\frac{d\mu_1}{d\lambda}}\right)=
W^{(\lambda)}\left(f_2\sqrt{\frac{d\mu_2}{d\lambda}}\right).
\end{equation}
Hence, for every $\varphi\in\mathbf L^2(\la)$ we have
\[
\begin{split}
\int_M\varphi(x)\left(
f_1(x)\sqrt{\frac{d\mu_1}{d\lambda}}(x)-f_2(x)
\sqrt{\frac{d\mu_2}{d\lambda}}(x)\right)d\lambda(x)&=\\
&\hspace{-5cm}= \mathbb
E_Q\left(W^{(\lambda)}(\varphi)\underbrace{W^{(\lambda)}
\left(f_1\sqrt{\frac{d\mu_1}{d\lambda}}-f_2
\sqrt{\frac{d\mu_2}{d\lambda}}d\lambda\right)}_{\mbox{$=0$}}\right)\\
&\hspace{-5cm}=0,
\end{split}
\]
as follows from \eqref{eq21}. Since this holds for all
$\varphi\in\mathbf L^2(\lambda)$ we conclude that
\[
f_1\sqrt{\frac{d\mu_1}{d\lambda}}
=f_2\sqrt{\frac{d\mu_2}{d\lambda}},\quad \lambda\,\, a.e.,
\]
that is $(f_1,\mu_1)\,\,\sim \,\, (f_2,\mu_2)$.
\mbox{}\qed\mbox{}\\

\section{The first main theorem}
\setcounter{equation}{0}
\label{sec4}

Starting with a measure space $M$ and a Borel sigma-algebra
$\mathscr B$, we get for every sigma-finite measure $\mu$  on $M$
an associated Gaussian process  $W^{(\mu)}$ . Now, for every
function  $f\in\mathbf L^2(\mu)$, we may therefore compute an
associated Ito-integral of $f$ with respect to this Gaussian
process  $W^{(\mu)}$; see Proposition \ref{rk3.2}. We denote this
Ito-integral by $W^{(\mu)}(f)$. We proved in Section \ref{sec3}
that, when $f$ and $\mu$ are given, then the Gaussian random
variable $W^{(\mu)}(f)$  depends only on the equivalence class of
the pair  $( f , \mu)$. As a result we are able to show (Theorem
\ref{tm:4.1}) that all the Gaussian processes  $W^{(\mu)}$  merge
together (via a sigma-Hilbert space) to yield a single Gaussian
Hilbert space in the sense of Definition \ref{def2.2}.

\begin{Dn}
\label{def5new}
Let $(M,\mathscr B)$ be fixed, and let $\mathscr H$ denote the
corresponding Hilbert space of sigma-functions; see Definition
\ref{def3.1}. For $\mu\in\mathscr(M,\mathscr B)$ we set
\begin{equation}
\label{eq51}
\mathcal H(\mu)=\left\{f\sqrt{d\mu}\,\,|\,\,
f\in\mathbf L^2(d\mu)\right\},
\end{equation}
and
\begin{equation}
\label{51new}
\mathcal H_1(\mu)=\left\{f\sqrt{d\mu}\,\,|\,\, f\in\mathbf
L^2(d\mu),\,\, |f|\le 1 \,\, \mu\,\,a.e.\right\},
\end{equation}
\end{Dn}

\begin{La}
Let $\mu\in\mathcal M(M,\mathscr B)$ be fixed. Then the
map
\begin{equation}
\label{opeT}
Tf=f\sqrt{d\mu}
\end{equation}
defines an isometrically isomorphic from $\mathbf L^2(M,\mu)$ onto
$\mathcal H(\mu)$
\end{La}

{\bf Proof:} It follows from Definition \ref{def3.1} that $T$ is
isometric. We claim that it is onto. Indeed, a pair $(g,\nu)$ is
in $\mathcal H(\mu)$ if and only if $\lambda=\mu+\nu$ satisfies
\begin{equation}
\label{eq5001}
g\sqrt{\frac{d\nu}{d\lambda}}=f\sqrt{\frac{d\mu}{d\lambda}},\quad
a.e.\,\,\lambda.
\end{equation}
We claim that
\begin{equation}
\label{eqTstar}
T^*(g,\nu)=f,
\end{equation}
where $f$ is as in \eqref{eq5001}. Indeed, for all
$\varphi\in\mathbf L^2(\mu)$ we have
\[
\begin{split}
\langle T\varphi,(g,\nu)\rangle_{\mathcal H(\mu)} &=
\int_M\varphi(x)g(x)\sqrt{\frac{d\nu}{d\lambda}\frac{d\mu}{d\lambda}}(x)d\lambda(x)\\
&=\int_M \varphi(x)f(x)\frac{d\mu}{d\lambda}(x)d\lambda(x)\\
&=\int_M \varphi(x)f(x)d\mu(x),
\end{split}
\]
and so $T^*(g,\nu)=f$ as claimed, and $TT^*={\rm Id}_{\mathcal
H(\mu)}$.
\mbox{}\qed\mbox{}\\

\begin{Tm}
\label{tm:4.1}
Let $\mathcal H$ be the sigma-Hilbert space of Definition
\ref{def3.1}. Let $(M,\mathscr B)$ be as in Section \ref{sec2}.
Then, the map
\begin{equation}
\label{eq22}
f\sqrt{d\mu}\quad\longrightarrow \quad W^{(\mu)}(f),
\end{equation}
defined for every $\mu\in\mathscr M(M,\mathscr B)$ and
$f\in\mathbf L^2(d\mu)$, extends to an isometry $F\mapsto
\widetilde{W}(F)$ from $\mathcal H$ into $\mathbf
L^2(\Omega_s,Q)$. Furthermore,
$\left\{\widetilde{W}(F)\right\}_{F\in\mathcal H}$ is a Gaussian
$\mathcal H$-process in the sense of Definition \ref{def2.2},
i.e.,
\begin{equation}
\label{eq24}
\mathbb E_Q  (\widetilde{W}(F_1)\widetilde{W}(F_2))=\langle F_1,
F_2\rangle_{\mathcal H},\quad \forall F_1,F_2\in\mathscr H.
\end{equation}
\end{Tm}

{\bf Proof:} Suppose that $F_i=f_i\sqrt{d\mu_i}$, $i=1,2$. Then, defining
\[
\widetilde W(F_i)=W^{(\mu_i)}(f_i),\quad i=1,2,
\]
(see \eqref{eq22}), identity \eqref{eq24} holds.  Indeed,
\[
\begin{split}
\mathbb E_Q\left(\widetilde{W}(F_1)\widetilde{W}(F_2)\right)&=
\mathbb E_Q\left({W^{(\mu_1)}}(f_1){W^{(\mu_2)}}(f_2)\right)\\
&=\mathbb
E_Q\left({W^{(\lambda)}}\left(f_1\sqrt{\frac{d\mu_1}{d\lambda}}\right)
{W^{(\lambda)}}
\left(f_2\sqrt{\frac{d\mu_2}{d\lambda}}\right)\right)\\
&=\int_M f_1(x)f_2(x)\sqrt{\frac{d\mu_1}{d\lambda}(x)\frac{d\mu_2}{d\lambda}(x)}
d\lambda(x)\\
&=\langle F_1,F_2\rangle_{\mathcal H},
\end{split}
\]
where, in the last step, we used \eqref{eq9} in the definition of the inner product in $\mathcal H$.\\

We now turn to the linearity of $\widetilde{W}$. For the sum in $\mathcal H$ we have equation \eqref{eq10}. Hence,
\[
\begin{split}
\widetilde{W}(F_1+F_2)&=W^{(\lambda)}\left(f_1\sqrt{\frac{d\mu_1}{d\lambda}}+f_2\sqrt{\frac{d\mu_2}{d\lambda}}\right)\\
&=W^{(\lambda)}\left(f_1\sqrt{\frac{d\mu_1}{d\lambda}}\right)+W^{(\lambda)}\left(f_2\sqrt{\frac{d\mu_2}{d\lambda}}\right)\\
&=\widetilde{W}(F_1)+\widetilde{W}(F_2).
\end{split}
\]
It remains to prove that $\widetilde{W}(\cdot)$ satisfies the
joint Gaussian property stated in Remark \ref{rk2.3}. We must
prove that if $F_i=f_i\sqrt{d\mu_i}$, $i=1,2,\ldots, n$, then the
joint distribution of
\[
(\widetilde{W}(F_1),\widetilde{W}(F_2),\ldots, \widetilde{W}(F_n))
\]
is the Gaussian random variable in $\mathbb R^n$ with zero mean
and covariance matrix $\left(\langle F_i,F_j\rangle_{\mathcal
H}\right)_{i,j=1}^n$. To see this, pick $\lambda\in\mathscr
M(M,\mathscr B)$ such that $\mu_1,\mu_2,\ldots, \mu_n$ are all
absolutely continuous with respect to $\lambda$ (for instance,
$\lambda=\sum_{i=1}^n\mu_i$). Then, in view of \eqref{eq9},
\begin{equation}
\label{eq26}
\langle F_i,F_j\rangle_{\mathcal H}=\int_M
f_i(x)f_j(x) \sqrt{\frac{d\mu_i}{d\lambda}
(x)\frac{d\mu_j}{d\lambda}(x)}d\lambda(x).
\end{equation}
But,
\[
\begin{split}
\mathbb E_Q\left(\widetilde{W}(F_i)\widetilde{W}(F_j)\right)&=\mathbb E_Q\left(W^{(\mu_i)}(f_i)W^{(\mu_j)}(f_j)\right)\\
&=
\mathbb E_Q\left(W^{(\lambda)}\left(f_i\sqrt{\frac{d\mu_i}{d\lambda}}\right)W^{(\lambda)}
\left(f_j\sqrt{\frac{d\mu_j}{d\lambda}}\right)\right)\\
&=\langle
f_i\sqrt{\frac{d\mu_i}{d\lambda}},f_j\sqrt{\frac{d\mu_j}{d\lambda}}\rangle_{\mathbf L^2(\lambda)},
\end{split}
\]
which is equal to the right hand-side of \eqref{eq26}, and leads
to the desired conclusion.\mbox{}\qed\mbox{}\\

We conclude this section with:

\begin{Pn}
$f\in\mathcal H_1(\mu)$ if and only if $f$ is the correlation
function for two copies of $W^{(\mu)}$.
\end{Pn}

{\bf Proof:} One direction is clear. Let $W_1^{(\mu)}$ and
$W_2^{(\mu)}$ be two correlated copies of $W^{(\mu)}$ and set
\[
\nu(A)=\mathbb E\left((W_1)_A^{(\mu)}(W_2)_A^{(\mu)}\right).
\]
$\nu$ is a signed measure, defining the correlation between the
two copies of $W^{(\mu)}$. By a Cauchy-Schwarz inequality, $\nu$
is absolutely continuous with respect to $\mu$. One then checks
that the Radon-Nikodym $f$ belongs to $\mathcal H_1(\mu)$.
Conversely, given $f\in\mathcal H_1(\mu)$, it suffices to define a
signed measure by
\[
\nu(A)=\int_A f(x)d\mu(x)
\]
to  define two correlated copies of $W^{(\mu)}$.

\mbox{}\qed\mbox{}\\

\section{Representation of $W^{(\mu)}$ in an arbitrary probability
space $(\Omega,\mathscr F, P)$}
\setcounter{equation}{0}
\label{sec6}
In this section we prove  that the infinite-product measure space (Theorem \ref{tm5.1}) is universal in the sense that
every measure space $(\Omega, \mathscr F, P)$  which carries some Gaussian
processes $W^{(\mu)}$, i.e., makes
$W^{(\mu)}$  into an $\mathbf L^2$ Gaussian process, can be computed directly
from the universal infinite-product measure
space. This is spelled out in Theorems \ref{tm6.1} in this section and in
Theorem \ref{tm7.1} in the next section.\\

In the previous section we have established a decomposition of the
Gaussian process $W^{(\mu)}$ as an expansion in a system of
i.i.d. $N(0,1)$ random variables. As before $(M,\mathscr B,\mu)$
is a given sigma-finite measure space.

\begin{Tm}
\label{tm6.1}
Let $W^{(\mu)}$ be represented in a probability space $\mathbf
L^2(\Omega,\mathscr F,\mathbb P)$, see Definition \ref{def2.1},
and let $\left\{\varphi_j\right\}_{j\in\mathbb N}$ be an
orthonormal basis in $\mathbf L^2(\mu)$. Then, there is a system
$\left\{Z_j\right\}_{j\in\mathbb N}$ of i.i.d. $N(0,1)$ random
variables such that
\begin{equation}
\label{eq6.1}
W_A^{(\mu)}=\sum_{j=1}^\infty\left(\int_A\varphi_j(t)d\mu(t)\right)Z_j(\cdot)
\end{equation}
holds almost everywhere on $\Omega$ with respect to $\mathbb P$.
\end{Tm}

{\bf Proof:} Assume that $A\in\mathscr B$ and $0<\mu(A)<\infty$.
We proceed in a number of steps.\\

STEP 1: {\sl The system $\left\{Z_j\right\}_{j\in\mathbb N}$
\begin{equation}
Z_j=W^{(\mu)}(\varphi_j)=\int_M\varphi_j(x)dW_x^{(\mu)}\quad\mbox{(as
an Ito integral)}
\label{eq38}
\end{equation}
is a family of i.i.d. $N(0,1)$ variables.}\\

To see this, we use the construction in Proposition \ref{rk3.2}. Indeed,
\[
\begin{split}
\mathbb E\left(Z_jZ_k\right)&=\mathbb E\left(W^{(\mu)}(\varphi_j)
W^{(\mu)}(\varphi_k)\right)\\
&=\langle\varphi_j,\varphi_k\rangle_{\mathbf L_2(\mu)}\\
&=\delta_{j,k},\quad \forall j,k\in\mathbb N,
\end{split}
\]
are the desired orthogonality condition. The rest of the assertion is clear.\\

STEP 2: {\sl We show that the sum on the right-hand-side of \eqref{eq38}
converges in the norm of
$\mathbf L^2(\Omega,\mathscr F,\mathbb P)$ and}

\[
\begin{split}
\mathbb E\left(|W_A^{(\mu)}-\sum_{j=1}^\infty\left(\int_A\varphi_j
(x)d\mu(x)\right)Z_j|^2\right)&=
\mathbb E\left(|W_A^{(\mu)}|^2\right)-\\
&\hspace{-2cm}-2\sum_{j=1}^\infty\left(\int_A\varphi_j(x)d\mu(x)
\right)\mathbb E\left(W_A^{(\mu)}Z_j\right)+\\
&\hspace{-2cm}
+\sum_{j=1}^\infty|\int_A\varphi_j(x)d\mu(x)|^2\\
&=\mu(A)-\sum_{j=1}^\infty|\int_A\varphi_j(x)d\mu(x)|^2\\
&=\mu(A)-\|\chi_A\|^2_{\mathbf L^2(\mu)}\\
&=\mu(A)-\mu(A)=0.
\end{split}
\]
\mbox{}\qed\mbox{}\\
\begin{Cy}
\label{cy6.2} Consider the space $(M,\mathscr B,\mu)$ as in the
previous theorem, and let $W^{(\mu)}$ be represented in some
probability space $\mathbf L^2(\Omega, \mathscr F,\mathbb P)$.
Then, some point $x_0\in M$ is an atom, i.e.
$\mu(\left\{x_0\right\})>0$, where $\left\{x_0\right\}$ denotes
the singleton, if and only if the ONB
$\left\{\varphi_j\right\}_{j\in\mathbb N}$in $\mathbf L^2(\mu)$
has $\varphi(x_0)$ well defined, the expansion \eqref{eq38}
contains a term
\begin{equation}
\label{eq41}
\mu(\left\{x_0\right\})\sum_{j=1}^\infty
\varphi(\left\{x_0\right\})Z_j,
\end{equation}
and
\begin{equation}
\label{eq42}
\sum_{j=1}^\infty(\varphi_j(x_0))^2=\frac{1}{\mu(\left\{x_0\right\})}
\end{equation}
\end{Cy}
{\bf Proof:} Functions $f\in\mathbf L^2(\mu)$ are determined only
point-wise a.e with respect to $\mu$, but if
$\mu(\left\{x_0\right\})>0$, the functions $f$ are necessarily
well defined at the point $x_0$, i.e., $f(x_0)$ is a uniquely
defined finite number. We apply this to the functions $\varphi_j$
in the $\mathbf L^2(\mu)$-ONB  from \eqref{eq38}. Hence, the
contributions to the two sides in \eqref{eq38} corresponding to
$A=\left\{x_0\right\}\in\mathscr B$ are as follows:
\begin{equation}
\label{eq43}
W_{\left\{x_0\right\}}^{(\mu)}(\cdot)=
\mu(\left\{x_0\right\})\sum_{j=1}^\infty
\varphi(\left\{x_0\right\})Z_j(\cdot).
\end{equation}
Taking norms in $\mathbf L^2(\Omega,\mathscr F,\mathbb P)$ we get
\[
\mu(\left\{x_0\right\})=\left(\mu(\left\{x_0\right\})\right)^2\sum_{j=1}^\infty
\left(\varphi(x_0)\right)^2,
\]
and the desired conclusion \eqref{eq42} follows.
\mbox{}\qed\mbox{}\\

\begin{Rk}
\label{rk6.3}
Some care must be exercised in asssigning the random variable
$W^{(\mu)}_A$ to sets $A\in\mathscr B$ with $\mu(A)=0$, or
$\mu(A)=\infty$: If $\mu(A)=0$, we may take $W^{(\mu)}_A$ to have
law the Dirac distribution $\delta_0$ on $\mathbb R$ at $x=0$. In
view of \eqref{eq6.1} one may alternatively set $W^{(\mu)}_A=0$
if $\mu(A)=0$. There are two conventions for dealing with the
random variable $X^{(\mu)}_A$ when $\mu(A)=\infty$. One involves a
renormalization, somewhat subtle. For other purposes, if
$\mu(A)=\infty$, we may simply take the random variable
$W_A^{(\mu)}$ to have the uniform distribution.
\end{Rk}
\mbox{}\\

We now explain the connections between the present construction
and the processes we built in \cite{MR2793121,ajnfao}.
\begin{Ap}
Let $\mathscr S$ denote the  Schwartz space  of smooth functions
on $\mathbb R$ with its Fr\'echet topology and let $\mu$ be a
Borel measure on $\mathbb R$ such that $\int_{\mathbb
R}\frac{d\mu(t)}{(1+t^2)^p}<\infty$ for some $p\in\mathbb N_0$
(where $\mathbb N_0=\mathbb N\cup\left\{0\right\}$). Then:\\
$(i)$ The function $F$:
\[
F(\varphi)= e^{-\frac{1}{2}\int_{\mathbb
R}|\widehat{\varphi}(t)|^2d\mu(t)}
\]
(where $\widehat{\varphi}$ denotes the Fourier  transform of
$\varphi$) is positive definite and continuous from $\mathscr S$
into $\mathbb R_+$, in the Fr\'echet topology. By Minlos' theorem
there exists a uniquely defined probability measure $\mathbb P$ on
the space $\mathscr S^\prime$ of tempered distributions such that
\[
F(\varphi)=\mathbb E_{\mathbb P}\left( e^{i\langle\cdot\, ,\,
\varphi\,\rangle}\right).
\]
$(ii)$ Furthermore we showed that there is a Gaussian process on
$\mathscr S^\prime$ with the Wiener measure such that
\[
F(\varphi)=\mathbb E_{\rm
Wiener}\left(e^{iX_\varphi^{(\mu)}}\right).
\]

From the results of the present paper, we then get
\[
F(\varphi)=\mathbb
E_Q\left(e^{iW^{(\mu)}(\widehat{\varphi})}\right),
\]
where $Q$ is the probability measure defined in Lemma
\ref{le5.1}, and where the process $W^{(\mu)}$ is constructed in
Proposition \ref{rk3.2} and Theorem \ref{tm5.1}.\\

In summary, the two Gaussian processes $X^{(\mu)}(\varphi)$ and
$W^{(\mu)}(\widehat{\varphi})$ have the same generating function.
\end{Ap}

As a corollary we have:

\begin{Cy}
Let $\lambda=dx$ denote the Lebesgue measure $dx$ on the real line
and the Gaussian processes $W^{(\mu)}(\varphi)$ constructed from
measures $\mu$ such that $\mu<<\lambda$ include the fractional
Brownian motion. We get this from the choice $d\mu(x)=c_H
|x|^{2H}dx$ where $H\in(0,1)$ and $c_H$ is some appropriate
constant.
\end{Cy}
For some recent work on the fractional Brownian motion, see also
\cite{AJP,al:pota,JSong1,Maslowski}.

\section{The probability space
$(\times_{\mathbb N} \mathbb R,\mathscr F,Q_K)$ is universal}
\setcounter{equation}{0}
\label{sec7}
Suppose a Gaussian
processes $W^{(\mu)}$  is represented in some measure space
$(\Omega, \mathscr F, P)$, we will then be able to compute the
measure $P$, and study how it depends on the initial measure
$\mu$ on $M$. This we do in
Theorem \ref{tm7.1} below, which also yields a measure-isomorphism
connecting $P$ to an infinite-product measure.\\

In this section we will show that when $(M,\mathscr B,\mu)$ is
given as above, that is, is some fixed sigma-finite measure
space, then every realization of the corresponding Gaussian
process $W^{(\mu)}$ factors through $(\times_{\mathbb N} \mathbb
R,Q_K)$. More precisely suppose that $W^{(\mu)}$ is realized as a
Gaussian $\mathbf L^2$-process in some probability space
$(\Omega,\mathscr F,\mathbb P)$, then there is a factor-mapping
setting up to an isomorphism of the respective Gaussian processes
on $(\Omega,\mathscr F,\mathbb P)$ and on $(\times_{\mathbb N}
\mathbb R,Q_K)$.

\begin{Tm}
\label{tm7.1}
Let $(M,\mathscr B,\mu)$ be fixed, and let the associated process
(see Definition \ref{def2.1} and Theorem \ref{tm5.1})
be realized in $\mathbf L^2(\Omega,\mathscr F,\mathbb P)$, where
$(\Omega,\mathscr F,\mathbb P)$ is a probability space.
Suppose
\begin{equation}
\label{eq44}
\mathscr F=\sigma-{\rm alg.}\left\{W_A^{(\mu)}\,\,|\,\, A\in\mathscr B\right\}.
\end{equation}
Then, the following assertions hold:\\
$(i)$ For all $A\in\mathscr B$ with $0<\mu(A)<\infty$ and $a,b\in\mathbb R$ (with $a<b$) we have
\begin{equation}
\label{eq45}
\mathbb P\left(\left\{\w\in\Omega\,\,|\,\, a<W_A^{(\mu)}(\w)\le b\right\}\right)=
\gamma_1\left((\frac{a}{\sqrt{\mu(A)}},\frac{b}{\sqrt{\mu(a)}}]\right),
\end{equation}
where $\gamma_1$ is the standard $N(0,1)$-Gaussian.\\
$(ii)$ There is a measure isomorphism
\[
\Psi\,\,\,:\,\,\, \Omega\longrightarrow\times_{\mathbb N}\mathbb R
\]
such that
\begin{equation}
\label{bgu:010812} \mathbb P\circ \Psi^{-1}=Q_K\quad{\rm and}\quad
W^{(Q_K,\mu)}\circ \Psi=W^{(\mu)}
\end{equation}
hold almost everywhere on $\Omega$, and where $W^{(Q,\mu)}$
denotes the realization of $W^{(\mu)}$ on $(\times_{\mathbb N}
\mathbb R,Q_K)$ from Section \ref{sec3}.
\end{Tm}

{\bf Proof:} Since $W_A^{(\mu)}\in\mathbf L^2(\Omega,\mathscr F,\mathbb P)$,
it follows that every cylinder set specified
as in $(i)$, i.e.,
\[
\left\{a<W_A^{(\mu)}\le b\right\},
\]
belongs to $\mathscr F$. since $W^{(\mu)}_A$ is a Gaussian
variable with law $N(0,\mu(A))$, formula \eqref{eq45}  from $(i)$
must hold. Now pick an ONB $\left\{\varphi_j\right\}_{j\in\mathbb
N}$ in $\mathbf L^2(\mu)$ and, following Theorem \ref{tm5.1}, set
\[
Z_j\, :=\, W^{(\mu)}(\varphi_j),\quad j\in\mathbb N.
\]
Then, $\left\{Z_j\right\}_{j\in\mathbb N}$ is an i.i.d. $N(0,1)$
family, and \eqref{eq38} holds. Now define
\[
\Psi\,\, :\,\, \Omega\quad \longrightarrow\quad \times_{\mathbb
N}\mathbb R
\]
by
\begin{equation}
\label{eq:psi}
\Psi(\w)=\left(Z_j(\w)\right)_{j\in\mathbb N},
\end{equation}
or equivalently,
\begin{equation}
\label{eq46}
\pi_j\circ \Psi=Z_j,\quad \forall j\in\mathbb N.
\end{equation}
Applying \eqref{eq46} to the expansion  \eqref{eq38} for
$W^{(\mu)}$ and for $W^{(Q_K,\mu)}$, we see get
\[
W^{(Q_K,\mu)}\circ \Psi=W^{(\mu)},
\]
that is, the stochastic process $W^{(\mu)}$ factors as stated.\\

Using again \eqref{eq38} from Theorem \ref{tm6.1}, we see that
\begin{equation}
\label{47} \Psi\left(\left\{a<W_A^{(\mu)}\le b\right\}\right)
\subseteq \left\{(\xi_j)_{j\in\mathbb N}\in\times_{\mathbb
N}\mathbb R\,\,\big|\,\,
a<\sum_{j=1}^\infty\xi_j\int_A\varphi_j(x)d\mu(x)\le b\right\},
\end{equation}
and that
\[
\begin{split}
Q\left(\Psi\left(\left\{a<W_A^{(\mu)}\le b\right\}\right)\right)&=
\gamma_1\left((\frac{a}{\sqrt{\mu(A)}},
\frac{b}{\sqrt{\mu(A)}}]\right)\\
&=\gamma_A\left((a,b]\right)\\
&=\mathbb
P\left(\left\{a<W_A^{(\mu)}\le b\right\}\right),
\end{split}
\]
where we used \eqref{eq45} in the last step of the reasoning.
Since $\mathscr F_s$ is generated (as a sigma-algebra)by the
cylinder sets, the final assertion $\mathbb P\circ \Psi^{-1}=Q_K$
in $(ii)$ follows. We do this by passing from monic subsets
$\left\{a<W_A^{(\mu)}\le b\right\}$, to finite functions, and to
measurable functions by inductive limit.\\

A function $F$ on $(\Omega,\mathscr F)$ is said to be finite if
there is $n\in\mathbb Z_+$, a bounded ${\mathbb R}^n$-Borel
function $f_n$, and $A_1,\ldots , A_n\in\mathscr B$ such that
\begin{equation}
\label{eq20111}
F(\cdot)=f_n(W_{A_1}^{(\mu)}(\cdot),\cdots,W_{A_n}^{(\mu)}(\cdot)).
\end{equation}
With $F$ as in \eqref{eq20111}, we then have
\[
\begin{split}
\mathbb E_{\mathbb
P}(F)&=\int_{\Omega}f_n(W_{A_1}^{(\mu)}(\w),\cdots,W_{A_n}^{(\mu)}(\w))d\mathbb
P(\w)\\
&=\underbrace{\iint\cdots\int}_{\mathbb R^n}f_n(x_1,x_2,\ldots,
x_n)\times\\
&\hspace{5mm}\times\gamma_{A_1}(x_1)\gamma_{A_2}(x_2-x_1)\cdots
\gamma_{A_n}(x_n-x_{n-1})dx_1dx_2\cdots dx_n\\
&=\underbrace{\iint\cdots\int}_{\mathbb
R^n}f_n(y_1,y_1+y_2,\ldots,
y_1+y_2+\cdots +y_n)\times\\
&\hspace{5mm}\times\gamma_{A_1}(y_1)\gamma_{A_2}(y_2)\cdots
\gamma_{A_n}(y_n)dy_1dy_2\cdots dy_n\\
&=\int_{\mathbf s^\prime}f_n(W_{A_1}^{(Q_K,\mu)},\cdots
W^{(Q_K,\mu)}_{A_n})dQ_K\\
&= \int_{\Omega}f_n(W_{A_1}^{(\mu)},\cdots
W^{(\mu)}_{A_n})d\left(Q\circ\Psi\right)\\
&=\mathbb E_{Q_K\circ \Psi}\left(F\right).
\end{split}
\]
\mbox{}\qed\mbox{}\\

\begin{Cy}
Fix $(M,\mathscr B,\mu)$ as in Theorems \ref{tm6.1} and \ref{tm7.1},
and denote by $\Omega_{\mathscr B}$
the set of all {\rm finitely} additive functions $\w\,\,:\,\,\mathscr B
\longrightarrow\mathbb R$. Set
\begin{equation}
\label{eq7.5}
W_A(\w)=\w(A),\quad \forall A\in\mathscr B,
\end{equation}
and let $\left\{Z_j\right\}_{j\in\mathbb N}$ be the corresponding i.i.d.
$N(0,1)$ system from Theorem \ref{tm6.1}
(The measure $\mathbb P$ on $\omega_{\mathscr B}$ is from \eqref{bgu:010812} in
Theorem \ref{tm7.1}, i.e.
\[
\mathbb P\circ \Psi^{-1}=Q_K,
\]
where $\Psi$ is given by \eqref{eq:psi}, that is,
$\Psi(\w)=\left(Z_j(\w)\right)_{j\in\mathbb N}$, $\forall
\w\in\mathscr B$). For $\xi=(\xi_j)_{j\in\mathbb N}\in\mathbf
s^\prime$, and $A\in\mathscr B$, set
\begin{equation}
\label{eq7.6}
\Gamma(\xi)(A)\,:=\,\sum_{j=1}^\infty\xi_j\int_A\varphi_j(x)d\mu(x).
\end{equation}
Then,
\begin{align*}
\hspace{-2cm}(i)\quad\hspace{2cm}&\Gamma(\xi)\in\Omega_{\mathscr B}&\\
\hspace{-2cm}(ii)\quad\hspace{2cm}&\Psi(\Gamma(\xi))=\xi,\quad \forall
\xi\in\mathbf s^\prime,&\\
\intertext{and} \hspace{-2cm}(iii)\quad\hspace{2cm}& \mathbb
P\left(\left\{\w\in\Omega_{\mathscr
B},\,\,|\,\,\Gamma\left(\Psi(\w)\right)=\w\right\}\right)=1&.
\end{align*}
\label{Cy7.2}
\end{Cy}

{\bf Proof:} The asserted conclusions follow from Theorem
\ref{tm6.1} and \ref{tm7.1}. Note that $(iii)$ in the Corollary
says that
\begin{equation}
\label{eq7.8}
\Gamma\circ \Psi={\rm Id}_{\Omega_{\mathscr
B}},\quad \mathbb P_{\mathscr B}\,\,\, a.e.
\end{equation}
where $\mathbb P_{\mathscr B}$ is the measure on
$\Omega_{\mathscr B}$ given by
\begin{equation}
\label{eq7.9}
\mathbb P_{\mathscr B}\circ \Psi^{-1}=Q_K.
\end{equation}

Now, formula \eqref{eq6.1} is an identity in $\mathbf
L^2(\Omega_{\mathscr B},\mathbb P_{\mathscr B})$. Since
\[
W^{(\mathscr B)}_A(\w)=\w(A),\quad\forall \w\in\Omega_{\mathscr
B},
\]
we get the following $\mathbb P_{\mathscr B}$-a.e. identity
holding on $\Omega_{\mathscr B}$:
\[
\begin{split}
\w(A)&=W^{(\mathscr B)}_A(\w)\\
&=\sum_{j=1}^\infty
\left(\int_A\varphi_j(x)d\mu(x)\right)Z_j(\w)\\
&= \left(\left(\Gamma\circ\Psi\right)(\w)\right)(A),\quad \forall
A\in\mathscr B.
\end{split}
\]
This proves \eqref{eq7.8}.
\mbox{}\qed\mbox{}\\

\begin{Cy}
\label{Cy7.31}
Let $\Psi\,\,:\,\,\Omega\quad\longrightarrow\quad \prod_{\mathbb
N}\mathbb R$ be as in \eqref{eq:psi}, and define the induced
operator $A$ from the  bounded Borel function defined on
$\prod_{\mathbb N}\mathbb R$ into the bounded Borel function
defined on $\Omega$,
\[
Af=f\circ \Psi.
\]
Then, $A$ is a Markov operator (see
\cite{arveson:markov}), i.e. the following properties hold:\\

$(i)$ $f\ge 0\,\, {\rm a.e.}\quad\Longrightarrow\quad Af\ge
0,\,\, {\rm a.e.}$,\\

$(ii)$ $A\mathbf 1=\mathbf 1$,\\

$(iii)$ $A^*\mathbf 1=\mathbf 1$.
\end{Cy}

{\bf Proof:} Note that in $(ii)$ and $(iii)$ the symbol $\mathbf
1$ denote the constant function equal to $1$ in the respective
measured spaces. Properties $(i)$ and $(ii)$ are clear. For
$(iii)$ we use the fact that functions of the form
\[
F(\xi)=f_n(\xi_1,\ldots, \xi_n).
\]
where $f_n$ is a bounded Borel function on $\mathbb R^n$ are
$\mathbf L^2$ dense. For  such a function, we want to check that
\[
\mathbb E_{Q_K}\left(F(\mathbf 1-A^*\mathbf 1)\right)=0,
\]
or, equivalently,
\[
\mathbb E_P\left(A(F)\right)=\mathbb E_{Q_T}(F).
\]
We have
\[
\begin{split}
\mathbb E_P\left(A(F)\right)&=\mathbb
E\left(f_n(Z_1(\cdot),\ldots,Z_n(\cdot)\right)\\
&= \underbrace{\iint\cdots\int}_{\mathbb R^n} f_n(x_1,\ldots,
x_n)\gamma_n(x_1,\ldots, x_n)dx_1\cdots dx_n\\
&=\mathbb E_{Q_K}(F),
\end{split}
\]
where we have used the properties from Theorem \ref{tm7.1} for
the respective measures $P$ and $Q_K$, as well as the i.i.d.
system $\left(Z_j\right)_{j\in\mathbb N}$ from \eqref{eq38} in
Theorem \ref{tm6.1}.
\mbox{}\qed\mbox{}\\
\begin{Cy}
\label{Cy7.3}
Let $(M,\mathscr B,\mu)$, $\Omega_{\mathscr B}$ and $\mathbb
P:=\mathbb P_{\mathscr B}$ be as in Corollary \ref{Cy7.2}. Set on
$\mathscr B\times \mathscr B$
\begin{equation}
\label{eq7.91}
K(A,B)=e^{\left(\mu(A\cap B)-\frac{\mu(A)+\mu(B)}{2}\right)}.
\end{equation}
Then we may define a Fourier transform $F\mapsto\widehat{F}$ from
$\mathbf L^2(\Omega_{\mathscr B},\mathbb P)$ onto the reproducing
kernel Hilbert space $\mathcal H(K)$ with reproducing kernel $K$
as in \eqref{eq7.91}. For $F\in\mathbf L^2(\Omega_{\mathscr
B},\mathbb P)$,
\begin{equation}
\label{eq710}
\widehat{F}(A)=\mathbb
E\left(F(\cdot)e^{iW_A^{(\mu)}(\cdot)}\right),\quad A\in\mathscr
B.
\end{equation}
Moreover the map $F\mapsto\widehat{F}$  is an isometric
isomorphism between the two Hilbert spaces.
\end{Cy}

{\bf Proof:} We begin  with finite sums of the form $\sum_{j\in J}
a_jK_{A_j}$, where the $a_j$ are real, $A_j\in\mathscr B$ and
$|J|<\infty$. Comparing the Hilbert norms we have
\[
\begin{split}
\| \sum_{j\in J} a_jK_{A_j}\|^2_{\mathcal H(K)}&=\sum_{j,k\in J}
a_ja_kK(A_j, A_k)\\
&=\sum_{j,k\in J}a_ja_k
e^{-\frac{1}{2}\|\chi_{A_j}-\chi_{A_k}\|^2_{\mathbf L^2(\mu)}}\\
&=\sum_{j,k\in J} a_ja_k\mathbb E\left(
e^{iW_{A_j}^{(\mu)}}e^{-iW_{A_k}^{(\mu)}}\right)\\
&=\|\sum_{j\in J}a_j e^{iW_{A_j}^{(\mu)}}\|^2_{\mathbf
L^2(\Omega_\mathscr B,\mathbb P)},
\end{split}
\]
where we have used Theorem \ref{tm5.1} (see in particular
\eqref{eq37}) in the last step in the computation. This complete
the proof of the isometry since such finite sums are dense in
$\mathcal H(K)$. Completing by taking $\mathcal H(K)$-norm
closure, we see that the adjoint of the map $J(F)=\widehat{F}$ is
isometric from $\mathcal H(K)$ into $\mathbf L^2(\Omega_\mathscr
B,\mathbb P)$. Indeed,
\begin{equation}
\label{eq711}
J\left(e^{iW_A^{(\mu)}}\right)=K_A\in\mathcal H(K),
\end{equation}
and
\begin{equation}
\label{712}
J^*(K_A)=e^{iW_A^{(\mu)}},\quad A\in\mathscr B.
\end{equation}
It remains to prove that
\[
\left\{e^{iW_A^{(\mu)}}\,\,|\,\,A\in\mathscr B\right\}
\]
span a dense subspace in $\mathbf L^2(\Omega_\mathscr B,\mathbb
P)$, and if $F\in\mathbf L^2(\Omega_{\mathscr B},\mathbb P)$ is
such that
\begin{equation}
\label{713}
\widehat{F}(A)=\mathbb
E\left(Fe^{iW^{(\mu)}_A}\right)=0,\quad\forall A\in\mathscr B,
\end{equation}
then $F=0$.\\

To verify this, we may use the known representation of $\mathbf
L^2(\Omega_{\mathscr B},\mathbb P)$ as the symmetric Fock space
over $\mathbf L^2(d\mu)$; see \cite{MR1244577}. We also make use
of Theorem \ref{tm7.1} above. Suppose $F\in\mathbf
L^2(\Omega_{\mathscr B},\mathbb P)$ satisfies \eqref{713}. In the
Fock-space representation,
\begin{equation}
\label{714}
F=\sum_{n=0}^\infty F_n
\end{equation}
is referring to Wiener chaos expansion of $F$, that is, the
orthogonal decomposition of $F$ along the orthogonal sum of all
symmetric $n$-tensors, as $n=0,1,2,\ldots$, and with $n=0$
referring to the vacuum vector. See also \cite{MR1978577}.
Substitution of \eqref{714} into \eqref{713} yields
\begin{equation}
\label{715}
\mathbb E\left(F_n\underbrace{W_A^{(\mu)}\times\cdots\times
W_A^{(\mu)}}_{\mbox{\rm $n$ times}}\right)=0,\quad\forall
A\in\mathscr B,\quad {\rm and}\quad n=0,1,\ldots
\end{equation}
Using now the Ito-integral from Proposition \ref{rk3.2}, equation
\eqref{715} may be rewritten as
\[
\underbrace{\int_A\int_A\cdots\int_A}_{\mbox{\rm $n$ times}}
F_n(x_1,x_2,\ldots , x_n)dW_{x_1}^{(\mu)} dW_{x_2}^{(\mu)}\cdots
dW_{x_n}^{(\mu)}=0,
\]
that is (and where $\otimes$ denotes the symmetric tensor
product),
\begin{equation}
\label{716}
F_n\quad \perp\quad \otimes_1^n\chi_A,\quad \forall A\in\mathscr
B.
\end{equation}
Since $F_n\in \otimes_1^n\mathbf L^2(\mu)$ is a symmetric tensor,
we conclude from \eqref{716} that $F_n=0$. This holds for
$n=0,1,\ldots$ and so by \eqref{714}, $F=0$.
\mbox{}\qed\mbox{}\\

\begin{Rk}
The fact that $K(A,B)$ is positive definite on $\mathscr B$ can be
checked also as follows: The function
\[
n(A,B)= -\mu(A\cap B)+\frac{\mu(A)+\mu(B)}{2}
\]
is conditionally negative on $\mathscr B$, and therefore the
function $e^{-n(A,B)}$ is positive definite there. See
\cite{MR86b:43001} for the latter.
\end{Rk}

\section{Iterated function systems}
\setcounter{equation}{0}
\label{sec8}
The purpose of the present section is to give an application of
the theorems from Sections \ref{sec6} and \ref{sec7} to iterated
function systems (IFS), see e.g. \cite{MR625600}. Such IFSs arise
in geometric measure theory, in harmonic analysis, and in the
study of dynamics of iterated substitutions with rational
functions (on Riemann surfaces); hence the name iterated function
system." With an IFS, we have the initial measure space $M$ and
a Borel sigma-algebra $\mathscr B$, coming with an additional
structure, a system of measurable endomorphisms. We will be
interested in those measures $\mu$ on $M$ which satisfy suitable
self-similarity properties with respect to the prescribed
endomorphisms in $M$. For background,
see e.g. \cite{MR1960424,MR2097020,MR2459320}.\\

Given a measure space $(M,\mathscr B)$ as in Section \ref{sec2},
i.e. $\mathscr B$ is a fixed Borel sigma-algebra of subsets of
$M$, by an {\sl iterated function system} (IFS), we mean a system
of endomorphisms $(\tau_i)_{i\in I}$
\[
\tau_i\,\,:\,\,  M\quad\longrightarrow\quad M,
\]
each $\tau_i$ assumed measurable and the index set $I$ usually finite.\\

If a family of measures $\mu$ on $\mathscr B$ is specified, each
$\tau_i$ is defined a.e.. Typically, $M$ will be a locally
compact Hausdorff space, and we assume that each $\tau_i$ is
continuous. The following restrictions will be placed on the
family $(\tau_i)_{i\in I}$:
\begin{equation}
\label{eq48}
\tau_i(M)\cap \tau_j(M)=\emptyset,\quad \forall i,j\in I\quad
(\mbox{\rm non-overlapping}),
\end{equation}
\begin{equation}
\bigcup_{i\in I} \tau_i(M)=M,\quad(\mbox{\rm cover}),
\label{eq49}
\end{equation}
and there is a measurable endomorphism $R$ from $M$ into $M$ such
that
\begin{equation}
\label{eq50}
R\circ \tau_i={\rm Id}_M,\quad \forall i\in I.
\end{equation}
We say that the family $\left\{\tau_i\right\}_{i\in I}$ is a
system of branches of an inverse to $R$. This is in particular
the case in applications to Riemann surfaces, where $R$ is
typically a rational function.\\

In view of the following definition, recall that we have defined
$\mathcal (\mu)$ in Definition \ref{def5new}.

\begin{Dn}
\label{def8.1}
Let $(M,\mathscr B)$ be fixed, and let $\mathscr H$ denote the
corresponding Hilbert space of sigma-functions; see Definition
\ref{def3.1}, and let $\mu\in\mathscr(M,\mathscr B)$.
%\begin{equation}
%\label{eq51} \mathcal H(\mu)=\left\{f\sqrt{d\mu}\,\,|\,\,
%f\in\mathbf L^2(d\mu)\right\}.
%\end{equation}
If $R\,\,:\,\, M\longrightarrow M$ is a measurable endomorphism,
we consider the measure $\mu\circ R^{-1}$, i.e.
\begin{equation}
\label{eq52}
(\mu\circ R^{-1})(A)=\mu(R^{-1}(A)),\quad \forall A\in\mathscr B,
\end{equation}
where
\[
R^{-1}(A)=\left\{x\in M\,\,|\,\, R(x)\in A\right\}.
\]
We set
\begin{equation}
\label{eq53}
\left( \mathcal H\circ R  \right)(\mu)=\left\{(f\circ
R)\sqrt{d\mu}\,\,|\,\, f\in \mathbf L^2(\mu\circ R^{-1})\right\}.
\end{equation}
\end{Dn}

\begin{Dn}
\label{def8.2}
Let $\mu\in\mathscr(M,\mathscr B)$, and let
$\left\{\tau_i\right\}_{i\in I}$ be as in
\eqref{eq48}-\eqref{eq50} in the previous definition. We say that
$(\mu,\left\{\tau_i\right\}_{i\in I})$ is an iterative function
system (IFS) if
\begin{equation}
\label{eq54}
\mu\circ \tau_i^{-1}<< \mu,\quad \forall i\in I.
\end{equation}
An IFS is said to be closed if
\begin{equation}
\label{eq55}
\sum_{i\in I}\frac{d(\mu\circ\tau_i^{-1})}{d\mu}=1.
\end{equation}
\end{Dn}

Note that the Radon-Nikodym derivatives in the summation
\eqref{eq45} are well defined on account of \eqref{eq54}.

\begin{Rk}
Special cases of $IFS$ have been widely studied in the
literature; see e.g. \cite{MR1960424,MR2097020,MR2459320,MR833073,MR1386842,MR1785620,MR2875207}.
\end{Rk}

In these examples, the Radon-Nikodym derivatives
$\frac{d(\mu\circ\tau_i^{-1})}{d\mu}$ in \eqref{eq55} are
constant functions, say
\[
\frac{d(\mu\circ\tau_i^{-1})}{d\mu}=p_i,\quad i\in I,
\]
and $\sum_{i\in I}p_i=1$, so that in particular $p_i\in(0,1)$. As
further special cases of this, we have the Cantor measures: For
example, let $M$ be the usual middle third Cantor set, and define
two endomorphisms
\[
\tau_0(x)=\frac{x}{3},\quad{\rm and}\quad \tau_1(x)=\frac{x+2}{3}.
\]
Then, there is a unique probability measure $\mu$ supported on
$M$ such that
\begin{equation}
\label{eq56}
\mu=\frac{1}{2}\left(\mu\circ\tau_0^{-1}+\mu\circ\tau_1^{-1}\right).
\end{equation}
This is an IFS, and $p_0=p_1=\frac{1}{2}$; compare with
\eqref{eq55}. The scaling dimension of $\mu$ is
$\log_3(2)=\frac{\ln 2}{\ln 3}$.

\begin{La}
\label{la6.4}
Let $(\mu,\left\{\tau_i\right\}_{i\in I})$ is an iterative
function system. Then for each $i\in I$ the mapping
\begin{equation}
\label{eq57}
(f,d\mu)\,\, \mapsto\,\, (f\circ R, \mu\circ \tau_i^{-1})
\end{equation}
induces (by passing to equivalence classes) an isometry from
$\mathcal H(\mu)$ into $(\mathcal H\circ
R)(\tau\circ\tau_i^{-1})$.
\end{La}

{\bf Proof:} In principle there are issues with passing the
transformation onto equivalence classes, but this can be dome via
an application of Lemma \ref{la3.4}. Hence in studying
\eqref{eq57}, the question reduces to checking instead that the
application
\begin{equation}
\label{eq58} W^{(\mu)}(f)\,\,\mapsto\,\, W^{(\mu\circ
\tau_i^{-1})}(f\circ R)
\end{equation}
is isometric. Indeed,
\[
\begin{split}
\mathbb E\left( |W^{(\mu)}(f)|^2\right)&=\int_M|f(x)|^2d\mu(x)\\
&=\int_M \left(|f\circ R\circ \tau_i(x)|^2\right)d\mu(x)\\
&=\int_M\left(|(f\circ R)(x)|^2\right)d(\mu\circ \tau_i^{-1})(x)\\
&=\mathbb E_Q\left(|W^{(\mu\circ\tau_i^{-1})}(f\circ R)|^2\right),
\end{split}
\]
which is the desired conclusion.
\mbox{}\qed\mbox{}\\

We now turn to representation of the Cuntz relations; see e.g.
\cite{BrJo02a,BrJo02b,MR2097020}.

\begin{Tm}
\label{tm8.5}
Let $(\mu,\left\{\tau_i\right\}_{i\in I})$ is a closed iterated
function system, and set
$g_i=\frac{d(\mu\circ\tau_i^{-1})}{d\mu}$ (see \eqref{eq54} and
\eqref{eq55}). Then the operators
\begin{equation}
\label{eq59}
S_i(f)=\chi_{\tau_i(M)}\sqrt{g_i}(f\circ R),\quad i\in I,
\end{equation}
define a representation of the Cuntz algebra $\mathcal O_I$ (with
index set $I$), acting on the Hilbert space $\mathbf L^2(\mu)$,
i.e. as isometries in $\mathbf L^2(\mu)$, the operators $S_i$ from
\eqref{eq59} satisfy:
\begin{eqnarray}
\label{eq601}
S_i^*S_j&=&\delta_{i,j} {\rm Id}_{\mathbf L^2(\mu)},\forall i,j\in I,\\
\sum_{i\in I} S_iS_i^*&=&{\rm Id}_{\mathbf L^2(\mu)}.
\label{eq602}
\end{eqnarray}
\end{Tm}

{\bf Proof:} Condition \eqref{eq601} is immediate from the
preceding lemma. Now fix $i\in I$. one checks that the $\mathbf
L^2(\mu)$-adjoint of the operator in \eqref{eq59} is
\begin{equation}
\label{eq61}
S_i^*\varphi=\varphi\circ\tau_i,\quad \forall \varphi\in \mathbf
L^2(\mu),\,\,\, \forall i\in I.
\end{equation}
We are now ready to verify \eqref{eq602}, i.e. the second Cuntz
relation. In this computation we make use of \eqref{eq55}, i.e.
\[
\sum_{i\in I} g_i=1,\quad \mu \,\,a.e.
\]
For $\varphi\in\mathbf L^2(\mu)$, we have:
\[
\begin{split}
\int_M|\varphi(x)|^2d\mu(x)&=\sum_{i\in
I}\int_M|\varphi(x)|^2g_i(x)d\mu(x)\\
&=\sum_{i\in I}\int_M|\varphi(x)|^2d(\mu\circ \tau_i^{-1})(x)\\
&=\sum_{i\in I}\int_M |\varphi\circ \tau_i|^2(x)d\mu(x)\\
&=\sum_{i\in I}\int_M |S_i^*\varphi|^2(x)d\mu(x)\\
&=\sum_{i\in I}\langle \varphi, S_iS_i^*\varphi\rangle_{\mathbf
L^2(\mu)}.
\end{split}
\]
Since this holds for all $\varphi\in\mathbf L^2(\mu)$ the desired
formula \eqref{eq602} has been verified.
\mbox{}\qed\mbox{}\\

\section{Gaussian versus non-Gaussian}
\label{sec9}
\setcounter{equation}{0}
In this section we show that
the theory, developed above, initially for Gaussian Hilbert
spaces, applies to some non-Gaussian cases; for example to those
arising in the study of random functions. To make this point
specific, we address such a problem for the special case of a
concrete random power series, studied as a family of infinite
Bernoulli convolutions on the real line. We know, see
\cite{MR0622034}, that every positive definite function may be
realized in a Gaussian Hilbert space. Our results in Sections
\ref{sec3}-\ref{sec5} are making this
precise in some settings dictated by applications to stochastic integration.\\

\begin{Dn}
\label{def9.1}
If $T$ is a set, then the function
\begin{equation}
\label{eq63}
C\,\, :\,\, T\times T\quad\longrightarrow\quad
\mathbb C
\end{equation}
is said to be {\sl positive semi-definite} (p.s.d) (we will also say {\sl positive definite}) if for
every finite subset $S\subset T$, and every family $\left\{a_s\right\}_{s\in S}\subset \mathbb C^{|S|}$, we have
\begin{equation}
\label{eq64}
\sum_{(s,t)\in S\times S} \overline{a_s}a_tC(s,t)\ge 0.
\end{equation}
\end{Dn}
A Gaussian representation of a p.s.d function consists of a
Hilbert space $\mathcal H$ and a function
\[
X\,\,:\,\, T\quad\longrightarrow\quad \mathcal H
\]
such that
\begin{equation}
\label{eq65}
C(s,t)=\langle X_s,X_t\rangle_{\mathcal H},\quad
\forall t,s\in T,
\end{equation}
such that, for all $t\in T$, $X_t$ is a Gaussian random variable
with zero mean, $\mathbb E(X_t)=0$, and moreover
\begin{equation}
\label{eq66}
\mathbb E\left(X_s^*X_t\right)=C(s,t).
\end{equation}
The following is an important example of a solution to the problem
\eqref{eq64}--\eqref{eq66}, when the Gaussian
restriction is relaxed. In its simplest form, it may be presented as follows:

\begin{Pn}
\label{lemma9.2}
Let $T=(0,1)$ and consider the function
\[
C\,\,:\,\, (0,1)\times(0,1)\quad\longrightarrow\quad \mathbb R^+
\]
defined by
\begin{equation}
\label{eq9.11}
C(\lambda,\rho)=\frac{\lambda\rho}{1-\lambda\rho}.
\end{equation}
There is a solution to the representation problem \eqref{eq65} in a
binary probability space $\Omega(2)=\times_{\mathbb N}
\left\{\pm1\right\}$ with the infinite coin-tossing probability product measure
\[
q:=\times_{\mathbb N}(\frac{1}{2},\frac{1}{2}).
\]
\end{Pn}

{\bf Proof:} We will be making use of facts on Bernoulli
convolutions. For some of the fundamentals in the theory of
Bernoulli convolutions, we refer to
\cite{MR2459320,MR1785620,MR1386842}. We consider on $\Omega(2)$
the system $\left\{\epsilon_k\right\}_{k\in\mathbb N}$ of random
variables
\[
\epsilon_k\left((\w_j)_{j\in\mathbb N}\right)=\w_k,\quad\forall k\in\mathbb N.
\]
Denoting the expectation with respect to $q$ by $\mathbb E_q(\cdot)$ we have
\begin{equation}
\label{eq68}
\mathbb E_q(\epsilon_k)=0,\quad{\rm and}\quad \mathbb E_q(\epsilon_j\epsilon_k)=\delta_{j,k},\quad\forall j,k\in\mathbb N.
\end{equation}
The system $\left\{\epsilon_k\right\}_{k\in\mathbb N}$ is therefore i.i.d., but non-Gaussian. For $\lambda\in(0,1)$, set
\begin{equation}
\label{eq69}
X_\lambda (\w)=\sum_{k=1}^\infty\epsilon_k(\w)\lambda^k,\quad\forall \w \in\Omega(2).
\end{equation}
Such an expression is called a random power series. Then the distribution
\begin{equation}
\label{eq70}
\mu_\lambda=q\circ X_\lambda^{-1},
\end{equation}
(i.e. $\mu_\lambda(A)=q(X_\lambda^{-1}(A))$ for all Borel subsets $A$ of $\Omega(2)$) is the infinite Bernoulli convolution
measure given by its Fourier transform
\begin{equation}
\label{eq71}
\widehat{\mu}(\xi)=\prod_{n=1}^\infty\cos\left(2\pi \lambda^n\xi\right),
\quad \forall\xi\in\mathbb R.
\end{equation}
Equivalently, if $\tau_\pm(x)=\lambda(x\pm 1)$, the $\mu_\lambda$ is the
unique measure defined on the Borel
sigma-algebra $\mathscr B$ of $\mathbb R$ by
\begin{equation}
\label{72}
\mu_\lambda=\frac{1}{2}\left(\mu_\lambda\circ
\tau_+^{-1}+\mu_\lambda\circ \tau_-^{-1}\right),
\end{equation}
see also \eqref{eq56}. Note that for every $\lambda\in(0,1)$,
$\mu_\lambda$ has compact support strictly contained
in the open interval $(-1,1)$. We now verify the covariance property
\begin{equation}
\label{eq73}
\mathbf E_q\left(X_\lambda X_\rho\right)=\frac{\lambda\rho}{1-\lambda\rho},
\quad\forall\lambda,\rho\in(0,1).
\end{equation}
In the left hand-side of \eqref{eq73} we substitute \eqref{eq69},
and we make use of the i.i.d. properties \eqref{eq68}. Then
\[
\mathbf E_q\left(X_\lambda X_\rho\right)=
\sum_{k=1}^\infty\lambda^k\rho^k=\frac{\lambda\rho}{1-\lambda\rho}.
\]
\mbox{}\qed\mbox{}\\

\begin{Tm}
\label{tm9.3} (Peres-Schlag-Solomyak and Peres-Solomyak,
\cite{MR1785620,MR1386842}) There is a Borel function
\[
D\,\, :\,\,
[\frac{1}{2},1)\times(-1,1)\quad\longrightarrow\quad\mathbb R^+
\]
such that the following properties hold for all $f\in
C_c\left([\frac{1}{2},1)\times(-1,1)\right)$:\\
$(i)$ The integral
\[
\iint_{[\frac{1}{2},1)\times(-1,1)}f(\lambda,
x)d\mu_\lambda(x)d\lambda
\]
is well defined, where $d\lambda$ denotes the standard Lebesgue
measure restricted to
$[\frac{1}{2},1)$,\\
and\\
$(ii)$ it holds that
\[
\iint_{[\frac{1}{2},1)\times(-1,1)}f(\lambda,
x)d\mu_\lambda(x)d\lambda=
\iint_{[\frac{1}{2},1)\times(-1,1)}f(\lambda, x)D(\lambda,
x)dxd\lambda.
\]
\end{Tm}

We first present some corollaries of this result.
\begin{Dn}
\label{dn9.6} Set
\begin{equation}
\label{eq901}
{\rm AC}_2=\left\{\lambda\in[\frac{1}{2},1)\,\,|\,\, \mbox{the
Radon-Nikodym derivative}\,\,\frac{d\mu_\lambda}{dx}\in\mathbf
L^2(dx)\right\}.
\end{equation}
(Note that the existence is part of the definition).
\end{Dn}

\begin{Rk}
\label{rk902}
The theorem asserts that ${\rm AC}_2$  has Lebesgue measure equal
to $1/2$, i.e. $\mu_\lambda$ is singular only on a subset of
$[\frac{1}{2},1)$ of measure zero. By a result of Erd\"os (see
\cite{MR0000311}), when $\lambda = g^{-1}$  where
$g=\frac{\sqrt{5}+1}{2}$ is the Golden ratio,, then $\mu_\lambda$
is singular. Otherwise it is absolutely continuous on a subset in
$[\frac{1}{2}, 1)$ of full measure.
\end{Rk}

\begin{Cy}
\label{cy903}
For a number $\lambda\in[\frac{1}{2},1)$, the following
conditions are equivalent:\\
$(i)$ The function
\[
t\mapsto \prod_{n=1}^\infty\cos(\lambda^nt),\quad t\in\mathbb R
\]
belongs to $\mathbf L^2(\mathbb R,dx)$.\\
$(ii)$ We have
\[
\liminf_{r\downarrow0}\frac{1}{2r}\mu_\lambda\left([x-r,x+x]\right)<\infty
\]
for a.a. $x\in\mathbb R$. In this  case, we may take
\begin{equation}
\label{eq904}
D(\lambda,
x)=\liminf_{r\downarrow0}\frac{1}{2r}\mu_\lambda\left([x-r,x+x]\right)\in
\mathbf L^2((-1,1),dx)
\end{equation}
in \eqref{eq901}.\\
\end{Cy}

\begin{Cy}
\label{cy904}
The points $\lambda\in{\rm AC}_2$ correspond to a single
equivalence class in the Hilbert space $\mathcal H$ of Definition
\ref{def3.1}. If $\lambda_1$ and $\lambda_2$ belong to ${\rm
AC}_2$,  we have
\begin{equation}
\label{eq906}
\langle f_1\sqrt{d\mu_{\lambda_1}},
f_2\sqrt{d\mu_{\lambda_2}}\rangle_{\mathcal
H}=\int_{-1}^1f_1(x)f_2(x) \sqrt{D(\lambda_1,x)D(\lambda_2,x)}dx,
\end{equation}
where $dx$ is the Lebesgue measure.
\end{Cy}

{\bf Proof:} Using \eqref{72}
\begin{equation}
\label{eq907}
\int\varphi(x)d\mu_\lambda(x)=\frac{1}{2}\left(
\int\varphi(\lambda(x+1))d\mu_\lambda(x)+
\int\varphi(\lambda(x-1))d\mu_\lambda(x)\right),
\end{equation}
and a recursive iteration leads to the representation
\begin{equation}
\label{eq908}
\begin{split}
\widehat{\mu_\lambda}(t)&=\int_{\mathbb
R}e^{-itx}d\mu_\lambda(x)\\
&=\mathbb E_q\left(e^{-itX_\lambda}\right)\\
&=\prod_{n=1}^\infty \cos(\lambda^nt),
\end{split}
\end{equation}
with the right hand-side of \eqref{eq908} converging point-wise
for all $t\in\mathbb R$.\\

If $\lambda\in{\rm AC}_2$, then
\[
D(\lambda, x)=\frac{d\mu_\lambda}{dx}\in\mathbf
L^2(-1,1)\subset\mathbf L^2(\mathbb R),
\]
and substitution into \eqref{eq908} yields
\[
\widehat{\mu_\lambda}(t)=\int_{\mathbb R}e^{-itx}D(\lambda,x)dx,
\]
and by the $\mathbf L^2(\mathbb R,dx)$-Fourier inversion,
\[
D(\lambda,x)=\int_{\mathbb
R}e^{itx}\prod_{n=1}^\infty\cos(\lambda^nt)dt
\]
for a.a. $x\in(-1,1)$. Hence, Plancherel's equality leads to
\[
\int_{-1}^1|D(\lambda,x)|^2dx=\int_{\mathbb
R}\prod_{n=1}^\infty\cos^2(\lambda^nt)dt<\infty.
\]
We now  turn to \eqref{eq906}. If $\lambda_1,\lambda_2\in{\rm
AC}_2$, then both $\mu_{\lambda_1}$ and $\mu_{\lambda_1}$ are
absolutely continuous with respect to Lebesgue measure, and by
\eqref{eq9} we get
\[
\begin{split}
\langle f_1\sqrt{d\mu_{\lambda_1}},
f_2\sqrt{d\mu_{\lambda_2}}\rangle_{\mathcal H}&=
\int_{-1}^1f_1(x)f_2(x)\sqrt{\frac{d\mu_{\lambda_1}}{dx}(x)
\frac{d\mu_{\lambda_2}}{dx}(x)}dx\\
& = \int_{-1}^1f_1(x)f_2(x) \sqrt{D(\lambda_1,x)D(\lambda_2,x)}dx.
\end{split}
\]
\mbox{}\qed\mbox{}\\

We showed that, when $\lambda$ is given in ${\rm AC}_2$, then the
corresponding  Bernoulli measure  $\mu_\lambda$ satisfies the
Bernoulli scaling law. But for  $\lambda$ fixed in ${\rm AC}_2$,
this then turns into a scaling identity for the $\mathbf L^2$
Radon-Nikodym derivative, a variant of the scaling law studied in
wavelet theory, but so far only for rational values of $\lambda$.
This fact is isolated in the corollary below. It is of interest
since there is very little known about $\mathbf L^2$ solutions to
scaling identity for non-rational values of $\lambda$. For the
literature on this we cite
\cite{BrJo02a,BrJo02b,MR1162107,MR1822853}.

\begin{Cy}
\label{cy10234}
Let $(\mu_\lambda)_{\lambda\in(0,1)}$ be the Bernoulli measures.
For $\lambda\in{\rm AC}_2$, let
$D(\lambda,\cdot)=\frac{d\mu_\lambda(x)}{dx}$ be the
Radon-Nikodym derivative. Extend $D(\lambda, x)$ to $x\in\mathbb
R$ by setting it to be equal to zero in the complement of
$(-1,1)$. Then,
\[
D_\lambda(\cdot)=D(\lambda,\cdot)\in\mathbf L_+^1(\mathbb R),\quad
\int_{\mathbb R}D(\lambda, x)dx=1,
\]
and
\begin{equation}
\label{eq12345}
D_\lambda(\lambda
x)=\frac{1}{2}\left(D_\lambda(x+1)+D_\lambda(x-1)\right)
\end{equation}
for a.a. $x$ with respect to the Lebesgue measure on $\mathbb R$.
\end{Cy}

{\bf Proof:} From the definition of ${\rm AC}_2$ we know that the
Radon-Nikodym derivative $x\mapsto D(\lambda, x)$ exists, and that
$D(\lambda,\cdot)\in\mathbf L_+^1(\mathbb R)\cap \mathbf
L^2(\mathbb R)$. Using \eqref{eq907}-\eqref{eq908} above, we
conclude that $\int_{\mathbb R} D(\lambda, x)dx=1$.
\mbox{}\qed\mbox{}\\

\begin{Rk}
\label{rk2001}
Note that for $\lambda=\frac{1}{2}$, equation \eqref{eq12345}
reduces to the standard scaling identity for the Haar wavelet
system in $\mathbf L^2(\mathbb R, dx)$. In wavelet theory, the
scaling identity is considered for $ N\in\mathbb Z_+$, $N>1$, as
follows: Given $N$, one studies solutions $\varphi\in\mathbf
L^2(\mathbb R,dx)$ to the scaling-rule
\[
\varphi(\frac{x}{N})=\sum_{k\in\mathbb Z}a_k\varphi(x-k),\quad
a.a. x,
\]
where $(a_k)_{k\in\mathbb Z}$ is square summable.
\end{Rk}

Before giving the proof of Theorem \ref{tm9.3} we need
preliminary lemmas:

\begin{La}
\label{la9.3}
\mbox{}\\
$(i)$ If $\lambda\in (0,\frac{1}{2})$, the measure $\mu_\lambda$
is singular with respect to Lebesgue measure, with scaling
dimension $D_s=-\frac{\ln2}{\ln\lambda}$, and the IFS defined by
$x\mapsto
\lambda(x\pm 1)$ is "non-overlapping".\\
$(ii)$ If $\lambda=\frac{1}{2}$, then $\mu_\lambda$ is equal to
the Lebesgue measure restricted to $[-1,1]$.\\
$(iii)$ For almost all $\lambda$ in $[\frac{1}{2},1)$, the measure
$d\mu_\lambda$ is absolutely continuous with respect to $dx$,
with Radon-Nikodym derivative
\[
\frac{d\mu_\lambda}{dx}(x)=D(\lambda, x)\in\mathbf L^2(-1,1).
\]
\end{La}

{\bf Proof:} The first two assertion are in the literature, and
$(iii)$ is from \cite{MR1386842}. It is our aim in Theorem
\ref{tm9.3} to give an independent proof in the reproducing kernel
\eqref{eq9.11} restricted to
$[\frac{1}{2},1)\times[\frac{1}{2},1)$; see also Proposition
\ref{lemma9.2} and equation \eqref{72}. \mbox{}\qed\mbox{}\\

Our purpose in connection with Theorem \ref{tm9.3} is as follows:
The proof of the result in \cite{MR1386842} relies on the
following estimate on $X_\lambda$  for a subset of points
$\la\in(\frac{1}{2},1)$, defined for measurable functions $F$ on
$\left(\times_{\mathbb N}\left\{\pm 1\right\}\right)\times
\left(\times_{\mathbb N}\left\{\pm 1\right\}\right)$, estimating
expectations
\begin{equation}
\label{eq12001}
\mathbb E_{q\times q}\left(\left(\mathbf 1\otimes
X_\la-X_\la\otimes\mathbf 1\right)F\right)
\end{equation}
where $\mathbf 1$ is the constant function $1$ on
$\times_{\mathbb N}\left\{\pm 1\right\}$.\\

One is in particular interested in \eqref{eq12001} in functions
$F$ of the form
\begin{equation}
\label{eqr}
F_r=\chi_{\left\{(\w,\w^\prime)\,\,\mbox{\tiny{\rm
such that}}\,\, |X_\la(\w)-X_\la(\w^\prime)|\le r\right\}},
\end{equation}
where $r\ge 0$.\\

For subintervals $J$ of $(\frac{1}{2},1)$ one must find estimate
on
\[
\int_J\mathbb E_{q\times q}( F_r)d\la
\]
In accomplishing this, the following three lemmas below are
helpful.

\begin{La}
\label{la9-11}
Let $\mathcal H$ be the reproducing kernel Hilbert
space from \eqref{eq9.11}, with $\lambda,\rho\in[\frac{1}{2},1)$,
and set
\[
k_\lambda(\rho)=\frac{\lambda\rho}{1-\lambda\rho}=
\langle
k_\lambda,k_\rho\rangle_{\mathcal H},\quad\forall
\lambda,\rho\in[\frac{1}{2},1).
\]
Then the assignment
\begin{equation}
\label{eq9.12}
k_\lambda\in\mathscr H\quad\mapsto\quad
X_\lambda(\cdot)\in \mathbf L^2(\times_{\mathbb
N}\left\{\pm1\right\},q)
\end{equation}
extends to a Hilbert space isometry of $\mathcal H$ into $\mathbf
L^2(\times_{\mathbb N}\left\{\pm1\right\},q)$.
\end{La}

{\bf Proof:} The conclusion follows from the basic axioms of reproducing
kernel Hilbert spaces once we verify that
\begin{equation}
\label{9.13}
\langle k_\lambda,k_\rho\rangle_{\mathcal H}=
\int_{\times_{\mathbb N}\left\{\pm
1\right\}}X_\lambda(\w)X_\rho(\w)dq(\w),
\end{equation}
equation \eqref{eq65} from the computation
\[
\langle k_\lambda,k_\rho\rangle_{\mathcal H}=\frac{\lambda\rho}{1-\lambda\rho}=
\mathbb E_q(X_\lambda X_\rho),
\]
by \eqref{eq9.11}.\mbox{}\qed\mbox{}\\

\begin{Dn}
\label{eeqhardy}
We denote by $\mathbf H^2(\mathbb D)$ the Hardy space of the
open unit disk of
functions analytic in the open unit disk $\mathbb D=\left\{z\in\mathbb C\,\,|\,\, |z|<1\right\}$ represented
as
\begin{equation}
f(z)=\sum_{n=0}^\infty a_nz^n,\quad \sum_{n=0}^\infty|a_n|^2<\infty,
\label{eq10012}
\end{equation}
and norm $\|f\|^2_{\mathbf H^2(\mathbb D)}=\sum_{n=0}^\infty|a_n|^2$, and set
\[
\mathbf H_0^2(\mathbb D)=\left\{f\in\mathbf H^2(\mathbb D)\,\,|\,\, f(0)=0\right\}.
\]
\end{Dn}

\begin{La}
\label{la9.5}
The reproducing kernel Hilbert space $\mathcal H$
from \eqref{eq9.11} is isometrically equal to $\mathbf
H_0^2(\mathbb D)$ via the map
\begin{equation}
\label{eq9.14}
k_\lambda\in\mathcal H\quad\mapsto\quad
\widetilde{k_\lambda}(z)=\sum_{n=1}^\infty \lambda^nz^n\in\mathbf H_0^2(\mathbb D).
\end{equation}
\end{La}
{\bf Proof:} It is immediate from the definition that the map
$k_\lambda\mapsto \widetilde{k_\lambda}$ in \eqref{eq9.14} extends to an isometry
\[
J\,\,:\,\,
\mathcal H\quad\longrightarrow\quad \mathbf H_0^2(\mathbb D).
\]
We claim that it is onto: ${\rm ran}~(J)=\mathbf H_0^2(\mathbb D)$.
Indeed, since $J$ is isometric, ${\rm ran}~J$ is
closed. now, if $f\in\mathbf H_0^2(\mathbb D)\ominus {\rm ran}~J$, then
\[
f(\lambda)=\langle f,\widetilde{k_\lambda}\rangle_{\mathbf H_0^2(\mathbb D)}=0,
\forall\lambda\in[\frac{1}{2},1).
\]
Since $f$ is analytic in $\mathbb D$ and $[\frac{1}{2},1)\subset \mathbb D$,
we conclude that $f\equiv 0$, and therefore
${\rm ran}~(J)=\mathbf H_0^2(\mathbb D)$ as claimed.
\mbox{}\qed\mbox{}\\

We now comment on the use of Lemmas \ref{la9-11} and \ref{la9.5}.
About \eqref{eq12001} the estimate
\[
\big|\mathbb E_{q\times q}\left(\left(\mathbf 1\otimes
X_\la-X_\la\otimes\mathbf
1\right)F\right)\big|\le\sqrt{\frac{2\la^2}{1-\la^2}}\|F\|_{\mathbf
L^2(q\times q)}
\]
follows form the Cauchy-Schwarz inequality, using that
\[
\mathbf 1\otimes X_\la\quad\perp\quad X_\la\otimes\mathbf 1
\]
in $\mathbf L^2(q\times q)$, and
\[
\|\mathbf 1\otimes X_\la\|^2_{\mathbf L^2(q\times q)}=
\|X_\la\|^2_{\mathbf L^2(q)}=\frac{\la^2}{1-\la^2}.
\]
See Proposition \ref{lemma9.2} and Lemma \ref{la9-11}.\\

As for estimating \eqref{eqr}, we make use of the Hardy space
representation in Lemma \ref{la9.5}. Under the isometry in
\eqref{9.13} the difference $|X_\la(\w)-X_\la(\w^\prime)|$ with
$\w_i=\w^\prime_i$ for $i=1,2,\ldots, k$ may be estimated in the
subspace $z^k\mathbf H_0^2(\mathbb D)$,, i.e. functions in
$\mathbf H^2(\mathbb D)$ vanishing at $0$ to order $k+1$.\\

\section{Boundaries of positive definite functions}
\setcounter{equation}{0} \label{sec10} In this section we apply
our results from Sections \ref{sec5} and \ref{sec7} into a general
boundary analysis for an arbitrarily given non-degenerate positive
definition function (Definition \ref{def9.1}). While it is known
that every non-degenerate positive definite function admits a
Gaussian representation, our construction here offers such a
representation in a form of a boundary in a sense which naturally
generalizes boundaries in classical analysis, for example
generalizing the known boundary analysis for the Szeg\"o kernel
of the disk. Again we stress that our starting point now is an
arbitrary fixed non-degenerate positive definite function $C$,
but $C$ is on $T\times T$ where $T$ may be any set, continuous or
discrete. For example $T$ may represent the vertices in some
infinite graph, and $C$ may be some associated energy form of the
graph $G$, induced by an electric network of $G$; see e.g.,
\cite{MR2643786,MR2799579}. A second recent application of
reproducing kernels and their RKHSs, is the theory of
(supervised) learning; see e.g.,
\cite{MR2944069,MR2810909,MR2488871}. The problem there is a
prediction of outputs based on observed samples; and for this the
kernel enters in representations of samples.
\\

Among the applications of stochastic processes, the theory of ``boundaries'' is noteworthy. Common to these is the need for
representations of functions on some set, say $T$, as integrals over some measure boundary space arising as a limiting
operation derived from the points in the initial set $T$. As example of this is the Hardy space $\mathbf H^2(\mathbb D)$
(see  Definition \ref{eeqhardy}), which is the reproducing kernel Hilbert space with kernel the Szeg\"o kernel
\begin{equation}
\label{1101}
C(z,w)=\frac{1}{1-zw^*},\quad z,w\in\mathbb D.
\end{equation}
If $\langle \cdot,\cdot\rangle_{\mathbf H^2(\mathbb D)}$ is the inner product of $\mathbf H^2(\mathbb D)$ we have
\begin{equation}
\label{1102}
f(w)=\langle f,C_w\rangle_{\mathbf H^2(\mathbb D)},\quad f\in\mathbf H^2(\mathbb D),\quad w\in\mathbb D.
\end{equation}
In this example we have
\begin{equation}
\label{eq10.4}
\frac{1}{1-zw^*}=\frac{1}{2\pi}\int_{-\pi}^\pi \frac{1}{1-ze^{-i\theta}}\frac{1}{1-w^*e^{i\theta}}d\theta.
\end{equation}
Now recall that the general case of positive definite functions $C$ on an arbitrary set $T$, as in Definition \ref{def9.1},
offers a generalization of the classical theory of the Hardy space recalled above. In this general case, the aim is
to provide a Gaussian measure space associated to an arbitrary given positive definite function
\begin{equation}
\label{eq10.5}
C\,\,:\,\, T\times T\quad\longrightarrow \quad\mathbb C.
\end{equation}
This measure space will be denoted by ${\rm bdr}_C(T)$, and it should be offer a direct
integral representation for \eqref{eq10.5} naturally generalizing \eqref{eq10.4}, where the boundary of $\mathbb D$ from
\eqref{eq10.4} is the circle $\left\{z\in\mathbb C\,\,|,\,\ |z|=1\right\}$.

\begin{Dn}
\label{def10.1} We say that a positive definite function $C$ on a
set $T$ is non-degenerate if the following two conditions
are satisfied:\\
$(i)$
\[
\quad{\rm dim}~\mathcal H(C)=\aleph_0,
\]
where $\mathcal H(C)$ is the reproducing kernel Hilbert space associated to $C$.\\
$(ii)$ The following implication holds:
\[
C(s,t_1)=C(s,t_2),\quad \forall s\in T\quad\Longrightarrow\quad t_1=t_2.
\]
\end{Dn}

\begin{Tm}
\label{tm10.2}
Let $C\,\,:\,\, T\times T\quad\longrightarrow\quad
\mathbb C$ be a non-degenerate positive definite function where
$T$ is some fixed set. Let $\mathbf s^\prime$ be the sequence
space introduced in Lemma \ref{le5.1} (see equation
\eqref{eq28}). Then there is a ${\rm weak}^*$-closed subspace
${\rm bdr}_C(T)\subset\mathbf s^\prime$, a Gaussian measure
$\mathbb P_C$ defined on the cylinder sigma-algebra in ${\rm
bdr}_C(T)$, and a Gaussian process $X$:
\begin{equation}
\label{eq10.6}
X_t\,\,:\,\, {\rm bdr}_C(T)\quad\longrightarrow\quad \mathbb C,\quad t\in T,
\end{equation}
such that $(i)$ we have
\begin{equation}
\label{eq10.7}
C(s,t)=\int_{{\rm bdr}_C(T)} X_s(\xi)^*X_t(\xi)d\mathbb P_C(\xi),\quad\forall s,t\in T,
\end{equation}
and, $(ii)$ $({\rm bdr}_C(T),\mathbb P_C,X_t)$ is a minimal solution to $(i)$.
\end{Tm}

{\bf Proof:} Let $\left\{\varphi_j\right\}_{j\in\mathbb N}$ be an orthonormal basis in $\mathcal H(C)$. It is well known that
\begin{equation}
\label{eq10.8}
C(t,s)=\sum_{j=1}^\infty \varphi_j(t)\varphi_j(s)^*,\quad\forall t,s\in T,
\end{equation}
and
\begin{equation}
\label{eq10.9}
\sum_{j=1}^\infty |\varphi_j(t)|^2=C(t,t)<\infty.
\end{equation}
Now define $\tau\,\, :\,\, T\,\,\longrightarrow\,\, \ell^2\subsetneq \mathbf s^\prime$ by
\begin{equation}
\label{eq10.10}
\tau(t)=\left(\varphi_j(t)\right)_{j\in\mathbb N},\quad t\in T.
\end{equation}
We claim that $\tau$ is one-to-one, and as result, we may identify points $t\in T$ with their image in
$\mathbf s^\prime$. Indeed, let $t_1,t_2\in T$ and suppose that $\tau(t_1)=\tau(t_2)$. Then,
\[
C(t,t_1)=\sum_{j=1}^\infty \varphi_j(t)(\varphi_j(t_1))^*=
\sum_{j=1}^\infty \varphi_j(t)(\varphi_j(t_2))^*=C(t,t_2),
\]
and in view of condition $(ii)$ in Definition \ref{def10.1} we conclude that $t_1=t_2$.\\

Set $\tau(T)=\left\{\tau(t)\,\,|\,\, t\in T\right\}$, and set ${\rm clo}_C(T)$ its closure in
$\mathbf s^\prime$. Here, by closure we mean the ${\rm weak}^*$-topology in $\mathbf s^\prime$ defined by the
duality between $\mathbf s$ and $\mathbf s^\prime$. The neighborhoods for this topology are generated by the cylinder
sets introduced in \eqref{eq30}. Finally, set
\begin{equation}
\label{eq10.11}
{\rm bdr}_C(T)={\rm clo}_C(T)\setminus\tau(T).
\end{equation}
Now, following Lemma \ref{le5.1}, set for $\xi\in{\rm bdr}_C(T)$
\begin{equation}
\label{eq10.12}
X_t(\xi)=\sum_{j=1}^\infty(\varphi_j(t))^*\pi_j(\xi)
=\sum_{j=1}^\infty\xi_j(\varphi_j(t))^*,
\end{equation}
the ``random'' function associated with the choice $\left\{\varphi_j\right\}$
of ONB in $\mathcal H(C)$.
Note that if $\xi$ in \eqref{eq10.12} is ``deterministic'', i.e.,
if there is a $s\in T$ such that
\[
\pi_j(\xi)=\xi_j=\varphi_j(s),\quad\forall j\in\mathbb N,
\]
then
\begin{equation}
\label{eq10.13}
X_t(\xi)=\sum_{j=1}^\infty \varphi_j(s)(\varphi_j(t))^*=C(t,s),\quad \forall t\in T.
\end{equation}

Now, define by $\mathbb P_C$ the measure on ${\rm bdr}_C(T)$ induced by $Q$ on $\mathbf s^\prime$, as in Theorem \ref{tm7.1}.
We get
\[
\begin{split}
\mathbb E_{\mathbb P_C}\left(X_t(\cdot)X_s(\cdot)^*\right)&=
\mathbb E_{\mathbb P_C}\left(\left(\sum_{j=1}^\infty \varphi_j(s)\pi_j^*\right)
\left(\sum_{k=1}^\infty \varphi_k(t)^*\pi_k\right)\right)\\
&=\sum_{j=1}^\infty\varphi_j(s)\varphi_j(t)^*\\
&=C(t,s),\quad\forall t,s\in T,
\end{split}
\]
whence the desired conclusion \eqref{eq10.7} in part $(i)$ of the theorem.
The other conclusion $(ii)$ follows form the assignment \eqref{eq10.11} in
the definition of ${\rm bdr}_C(T)$.
\mbox{}\qed\mbox{}\\

\begin{Ap}
Our boundary construction applies to electrical networks as
follows (see \cite{MR2862151}).\\

An electrical network is an infinite graph $(V, E, c)$ , $V$ for
vertices,  and $E$ for edges, where $c$ is a positive function on
$E$ , representing conductance. As sketched in \cite{MR2862151},
we get a reproducing kernel Hilbert space from the energy form
of  $(V, E, c)$. In \cite{MR2862151}, the authors propose one
boundary construction, and one can verify that the one from our
present Theorem \ref{tm10.2} applied to $\mathcal H$ is a
refinement.
\end{Ap}

\begin{Rk}
\label{rk10.3} Our construction of ${\rm bdr}_C(T)$ depends on
the choice of ONB in \eqref{eq10.8}, but the arguments in the
proof in Theorem \ref{tm10.2} above) show that two choices of ONB
$\left\{\varphi_j\right\}_{j\in\mathbb N}$ and
$\left\{\psi_k\right\}_{k\in\mathbb N}$ yields the same ${\rm
bdr}_C(T)$ if and only if there is an infinite unitary matrix
$(U_{j,k})_{(j,k)\in\mathbb N^2}$ such that
\[
\begin{split}
(i)& \hspace{1cm}\varphi_j=\sum_{k\in\mathbb N}U_{j,k}\psi_k,\\
\intertext{and the following equivalence holds:} (ii)
&\hspace{1cm} (b_j)_{j\in\mathbb N}\in\mathbf s\quad\iff\quad
(c_j)_{j\in\mathbb N}\in\mathbf s,\,\,{with}\,\,
c_j=\sum_{k\in\mathbb N} U_{j,k}b_k.
\end{split}
\]
In other words, the matrix-operation defined from $U$ preserves the
sequence space $\mathbf s$ of \eqref{27}.
\end{Rk}

\begin{Ex}
Consider now the Hardy space $\mathbf H^2(\mathbb D)$ (see Definition \ref{eeqhardy}).
On may check that, with the choice of the standard ONB in $\mathbf H^2(\mathbb D)$
\[
\varphi_k(z)=z^k,\quad k\in\mathbb
N_0:=\left\{0\right\}\cup\mathbb N,\quad z\in\mathbb D,
\]
we get
\[
{\rm bdr}_{\mbox{\rm
Szeg\"o}}=\left\{\left(e^{ik\theta}\right)_{k\in\mathbb
N_0}\,\,|\,\, \theta\in(-\pi,\pi]\right\},
\]
which by identification yields $(\pi, \pi]$, which is consistent
with \eqref{eq10.4} above.
\label{ex10.4}
\end{Ex}

\begin{Ex} (see \cite{aj-cras}).
\label{exemple:iowa1}
Here the pair $(C,T)$ from Definition \ref{def10.1} is as
follows: Consider the rational function $R(z)$ given by
\[
R(z)=z^4-2z^2,\quad z\in\mathbb C.
\]
Set $R_0(z)=z$, $R_1(z)=R(z)$ and
\[
R_n(z)=\underbrace{(R\circ R\circ \cdots\circ R)}_{\mbox{\rm $n$
times}}(z).
\]
Now set
\[
T=\Omega=\left\{z\in\mathbb C\,\,\mbox{{\rm such that}}\,\,
(R_n(z))_{n\in\mathbb N_0}\in\ell^1\right\},
\]
(where $\mathbb N_0=\mathbb N\cup\left\{0\right\}$) and on
$\Omega\times\Omega$ set
\[
C(z,w)=\prod_{n=0}^\infty(1+R_n(z)R_n(w)^*).
\]
\end{Ex}

Using the ideas of Exemple \ref{ex10.4} and from \cite{aj-cras}
we note that for this $(C,T)$ we get that ${\rm clo}_C(T)$ is the
filled Julia set of $R$. See also\ cite{MR1128089} for basic
properties of Julia sets derived from fixed rational functions of
a single complex variable.

\begin{Dn}
\label{dn10.5}
Let $(C,T)$ be as in Definition \ref{def10.1}. Following
\cite{MR2488871}, we say that $C$ is a {\rm Mercer-kernel} if:\\
$(i)$ $T$ is a compact metric space (with respect to some metric,
say $d$), and\\
$(ii)$ The function
\[
C\,\,:\,\, T\times T\quad\longrightarrow\quad \mathbb C
\]
is continuous with respect to $d\times d$.
\end{Dn}

\begin{Pn}
\label{pn10.6}
If $(C,T)$ is a Mercer kernel, then ${\rm clo}_C(T)=\tau(T)$; in
other words $\tau(T)$ from \eqref{eq10.10} is closed.
\end{Pn}

{\bf Proof:} Let $\xi\in\mathbf s^\prime$, and let
$(t_k)_{k\in\mathbb N}$ be a sequence of points of $T$ such that
$\lim_{k\rightarrow\infty} \tau(t_k)=\xi$; see the discussion
before Lemma \ref{le5.1}. Using $(i)$ in Definition \ref{dn10.5},
we may, without loss of generality, assume that the sequence
$(t_k)_{k\in\mathbb N}$  is convergent in $T$, i.e.
$\lim_{k\rightarrow\infty} d(t,t_k)=0$ where $t\in T$ is its
limit point.\\

Let $\left\{\varphi_k\right\}_{j\in\mathbb N}$ be an ONB in
$\mathcal H(C)$, see \eqref{eq10.8} in the proof of Theorem
\ref{tm10.2}. Then,
\[
\begin{split}
\|\tau(t)-\tau(t_k)\|^2_2&=\sum_{j=1}^\infty |\varphi_j(t)
-\varphi_j(t_k)|^2\\
&=C(t,t)-2{\rm Re}~C(t_k,t)+C(t_k,t_k),
\end{split}
\]
where we have used \eqref{eq10.8}-\eqref{eq10.10} in this
computation.\\

By virtue of Condition $(ii)$ in Definition \ref{dn10.5}, we now
note that the right hand-side in the last term converges to zero
as $k\rightarrow\infty$. But convergence in $\ell^2$ of the
sequence $\left(\tau(t_k)\right)_{k\in\mathbb N}$ implies
convergence in $\mathbf s^\prime$. We conclude that
$\tau(t)=\xi$, and so $\tau(T)$ is closed in $\mathbf s^\prime$.
\mbox{}\qed\mbox{}\\

\begin{Ex}
\label{exempleBM}
Let $T=I=[0,1]$ be the closed unit interval, and set
\[
C(t,s)=t\wedge s,\quad t,s\in I.
\]
Set
\begin{equation}
\label{varphi}
\varphi_k(t):\begin{cases} \sqrt{2}\,\frac{\sin
k\pi t}{k\pi},\quad k\in\mathbb N,\\
\,\,t,\hspace{1.3cm}\quad k=0.
\end{cases}
\end{equation}
Then:\\
$(i)$ $\tau(t):=\left(\varphi_k(t)\right)_{k\in\mathbb N_0}$
satisfies
\begin{equation}
\label{dist}
\|\tau(t)-\tau(s)\|_2^2=|t-s|,\quad t,s\in [0,1].
\end{equation}
$(ii)$ The map $t\mapsto \tau(t)$ is an  homeomorphism from $I$
onto a closed curve starting at $v_0=(0,0,0,\ldots)$ and with
endpoint
$v_1=(1,0,0,\ldots)$ in $\ell^2$.\\
$(iii)$ The curve in $(ii)$ has no self-intersection.
\end{Ex}

{\bf Proof of the claims in Exemple \ref{exempleBM}:} The
conclusions are  immediate from Proposition \ref{pn10.6}. Indeed,
the reproducing kernel Hilbert space associated to $C$ is
\begin{equation}
\label{sobolev}
\mathcal H=\left\{f\in\mathbf L^2(I)\,\,\big|\,\,
f^\prime\in\mathbf L^2(I)\,\,{\rm and}\,\, f(0)=0\right\},
\end{equation}
and one easily checks that the function system
$\left(\varphi_k\right)_{k\in\mathbb N_0}$ is an orthonormal
basis in $\mathcal H$. Indeed, for $j,k\in\mathbb N$,
\[
\langle\varphi_j,\varphi_k\rangle_{\mathcal H}=2\int_0^1\cos
(j\pi x)\cos(k\pi x)dx=\delta_{j,k}.
\]
The assertions follow then from Proposition \ref{pn10.6}. In this
example the Gaussian process from \eqref{eq10.12} associated with
$(C,I)$ is the Brownian motion. Hence
\begin{equation}
\label{eqBM}
\|\tau(t)-\tau(s)\|_2^2=\mathbb
E\left(|X_t-X_s|^2\right)=|t-s|,\quad t,s\in I,
\end{equation}
which is $(i)$, and also leads to $(ii)$ since $\tau$ is
one-to-one and continuous between two compact spaces, and so is an
homeomorphism. To justify \eqref{eqBM} note that the Hilbert norm
in $\mathcal H$ is $\|f\|^2_{\mathcal
H}=\int_0^1|f^\prime(x)|^2dx$.\\

Setting
\[
C_t(x)=\begin{cases} 0,\quad x<0,\\
                     x,\quad 0\le x\le t,\\
                     t,\quad t<x,
                     \end{cases}
\]
we get
\[
\begin{split}
\langle C_t,C_s\rangle_{\mathcal
H}&=\int_0^1\chi_{[0,t]}(x)\chi_{[0,s]}(x)dx\\
&=t\wedge s\\
&=C(t,s),
\end{split}
\]
and
\[
C(t,s)=ts+\frac{2}{\pi^2}\sum_{k=1}^\infty \frac{\sin (k\pi
t)\sin(k\pi s)}{k^2}.
\]
Finally, if there exist $t_1$ and $t_2$ in $(0,1)$ such that
$\tau(t_1)=\tau(t_2)$, then $C(t,t_1)=C(t,t_2)$  which is not
possible unless $t_1=t_2$.
\mbox{}\qed\mbox{}\\

{\bf Acknowledgments:} D. Alpay thanks the Earl Katz family for
endowing the chair which supported his research. The research of
the authors was supported in part by the Binational Science
Foundation grant 2010117, and Palle Jorgensen thanks the
department of mathematics for hospitality. We wish also to thank
our colleagues Dorin Dutkay (U.C. Florida), Rob Martin (Cape Town
University), Paul Muhly (University of Iowa), Judy Packer
(University of Colorado), Steen Pedersen (Wright St Univ),
Myung-Sin Song (University of South Illinois), Feng Tian (Wright
St University) for discussions.
\bibliographystyle{plain}
%\bibliography{/users/faculty/math/dany/bib/all}
%\bibliography{all}

\begin{thebibliography}{10}

\bibitem{AJP}
S.~Albeverio, P.E.T. Jorgensen, and A.M. Paolucci.
\newblock On fractional {B}rownian motion and wavelets.
\newblock {\em Complex analysis and Operator Theory}, 6:33--63, 2012.

\bibitem{MR2002b:47144}
D.~Alpay.
\newblock {\em The {S}chur algorithm, reproducing kernel spaces and system
  theory}.
\newblock American Mathematical Society, Providence, RI, 2001.
\newblock Translated from the 1998 French original by Stephen S. Wilson,
  Panoramas et Synth\`eses.

\bibitem{aal2}
D.~Alpay, H.~Attia, and D.~Levanony.
\newblock On the characteristics of a class of {G}aussian processes within the
  white noise space setting.
\newblock {\em Stochastic processes and applications}, 120:1074--1104, 2010.

\bibitem{aj-cras}
D.~Alpay and P.~Jorgensen.
\newblock {Espaces \`a noyau reproduisant de fonctions analytiques pour les
  ensembles de Julia remplis}.
\newblock Preprint. 2012.

\bibitem{ajnfao}
D.~Alpay and P.~Jorgensen.
\newblock Stochastic procesees induced by singular operators.
\newblock {\em {Numerical Functional Analysis and Optimization}}, 33:708--735,
  2012.

\bibitem{MR2793121}
D.~Alpay, P.~Jorgensen, and D.~Levanony.
\newblock A class of {G}aussian processes with fractional spectral measures.
\newblock {\em J. Funct. Anal.}, 261(2):507--541, 2011.

\bibitem{ajlm2}
D.~Alpay, P.~Jorgensen, I.~Lewkowicz, and I.~Marziano.
\newblock {Infinite product representations for kernels and iterations of
  functions}.
\newblock Preprint. 2012.

\bibitem{al:pota}
D.~Alpay and D.~Levanony.
\newblock On the reproducing kernel hilbert spaces associated with the
  fractional and bi--fractional brownian motions.
\newblock {\em {Potential Analysis}}, 28:163--184, 2008.

\bibitem{alp}
D.~Alpay, D.~Levanony, and A.~Pinhas.
\newblock Linear stochastic state space theory in the white noise space
  setting.
\newblock {\em {SIAM} {Journal of Control and Optimization}}, 48:5009--5027,
  2010.

\bibitem{arveson:markov}
W.~Arveson.
\newblock {Markov operators and $OS$-positive processes}.
\newblock {\em {Journal of Functional Analyis}}, 66:173--234, 1986.

\bibitem{MR1978577}
William Arveson.
\newblock {\em Noncommutative dynamics and {$E$}-semigroups}.
\newblock Springer Monographs in Mathematics. Springer-Verlag, New York, 2003.

\bibitem{MR86b:43001}
C.~Berg, J.~Christensen, and P.~Ressel.
\newblock {\em Harmonic analysis on semigroups}, volume 100 of {\em Graduate
  Texts in Mathematics}.
\newblock Springer-Verlag, New York, 1984.
\newblock Theory of positive definite and related functions.

\bibitem{BrJo02a}
O.~Bratteli and P.~Jorgensen.
\newblock {\em Wavelets through a looking glass}.
\newblock Applied and Numerical Harmonic Analysis. Birkh\"auser Boston Inc.,
  Boston, MA, 2002.

\bibitem{BrJo02b}
Ola Bratteli and Palle E.~T. Jorgensen.
\newblock Wavelet filters and infinite-dimensional unitary groups.
\newblock In {\em Wavelet analysis and applications ({G}uangzhou, 1999)},
  volume~25 of {\em AMS/IP Stud. Adv. Math.}, pages 35--65. Amer. Math. Soc.,
  Providence, RI, 2002.

\bibitem{MR1162107}
Ingrid Daubechies.
\newblock {\em Ten lectures on wavelets}, volume~61 of {\em CBMS-NSF Regional
  Conference Series in Applied Mathematics}.
\newblock Society for Industrial and Applied Mathematics (SIAM), Philadelphia,
  PA, 1992.

\bibitem{MR2133804}
Dorin~Ervin Dutkay and Palle E.~T. Jorgensen.
\newblock Hilbert spaces of martingales supporting certain
  substitution-dynamical systems.
\newblock {\em Conform. Geom. Dyn.}, 9:24--45 (electronic), 2005.

\bibitem{MR2643786}
Dorin~Ervin Dutkay and Palle E.~T. Jorgensen.
\newblock Spectral theory for discrete {L}aplacians.
\newblock {\em Complex Anal. Oper. Theory}, 4(1):1--38, 2010.

\bibitem{MR0000311}
Paul Erd{\"o}s.
\newblock On a family of symmetric {B}ernoulli convolutions.
\newblock {\em Amer. J. Math.}, 61:974--976, 1939.

\bibitem{MR0265548}
Leonard Gross.
\newblock Abstract {W}iener measure and infinite dimensional potential theory.
\newblock In {\em Lectures in {M}odern {A}nalysis and {A}pplications, {II}},
  pages 84--116. Lecture Notes in Mathematics, Vol. 140. Springer, Berlin,
  1970.

\bibitem{MR1244577}
T.~Hida, H.~Kuo, J.~Potthoff, and L.~Streit.
\newblock {\em White noise}, volume 253 of {\em Mathematics and its
  Applications}.
\newblock Kluwer Academic Publishers Group, Dordrecht, 1993.
\newblock An infinite-dimensional calculus.

\bibitem{MR2444857}
T.~Hida and Si~Si.
\newblock {\em Lectures on white noise functionals}.
\newblock World Scientific Publishing Co. Pte. Ltd., Hackensack, NJ, 2008.

\bibitem{Hida_BM}
Takeyuki Hida.
\newblock {\em Brownian motion}, volume~11 of {\em Applications of
  Mathematics}.
\newblock Springer-Verlag, New York, 1980.
\newblock Translated from the Japanese by the author and T. P. Speed.

\bibitem{MR2083706}
Takeyuki Hida.
\newblock A frontier of white noise analysis, in line with {I}t\^o calculus.
\newblock In {\em Stochastic analysis and related topics in {K}yoto}, volume~41
  of {\em Adv. Stud. Pure Math.}, pages 111--119. Math. Soc. Japan, Tokyo,
  2004.

\bibitem{MR1203453}
Takeyuki Hida, Hui-Hsiung Kuo, and Nobuaki Obata.
\newblock Transformations for white noise functionals.
\newblock {\em J. Funct. Anal.}, 111(2):259--277, 1993.

\bibitem{MR1408433}
H.~Holden, B.~{\O}ksendal, J.~Ub{\o}e, and T.~Zhang.
\newblock {\em Stochastic partial differential equations}.
\newblock Probability and its Applications. Birkh\"auser Boston Inc., Boston,
  MA, 1996.

\bibitem{MR625600}
John~E. Hutchinson.
\newblock Fractals and self-similarity.
\newblock {\em Indiana Univ. Math. J.}, 30(5):713--747, 1981.

\bibitem{MR1960424}
P.~E.~T. Jorgensen and D.~W. Kribs.
\newblock Wavelet representations and {F}ock space on positive matrices.
\newblock {\em J. Funct. Anal.}, 197(2):526--559, 2003.

\bibitem{MR2097020}
Palle E.~T. Jorgensen.
\newblock Iterated function systems, representations, and {H}ilbert space.
\newblock {\em Internat. J. Math.}, 15(8):813--832, 2004.

\bibitem{MR2459320}
Palle E.~T. Jorgensen, Keri Kornelson, and Karen Shuman.
\newblock Orthogonal exponentials for {B}ernoulli iterated function systems.
\newblock In {\em Representations, wavelets, and frames}, Appl. Numer. Harmon.
  Anal., pages 217--237. Birkh\"auser Boston, Boston, MA, 2008.

\bibitem{MR2862151}
Palle E.~T. Jorgensen and Erin P.~J. Pearse.
\newblock Gel\cprime fand triples and boundaries of infinite networks.
\newblock {\em New York J. Math.}, 17:745--781, 2011.

\bibitem{MR2799579}
Palle E.~T. Jorgensen and Erin P.~J. Pearse.
\newblock Spectral reciprocity and matrix representations of unbounded
  operators.
\newblock {\em J. Funct. Anal.}, 261(3):749--776, 2011.

\bibitem{JSong1}
{Jorgensen, Palle E. T. and Song, Myung-Sin}.
\newblock An extension of {W}iener integration with the use of operator theory.
\newblock {\em {J. Math. Phys.}}, 50(10):{103502, 11}, 2009.

\bibitem{MR833073}
Jean-Pierre Kahane.
\newblock {\em Some random series of functions}, volume~5 of {\em Cambridge
  Studies in Advanced Mathematics}.
\newblock Cambridge University Press, Cambridge, second edition, 1985.

\bibitem{Ka48}
Shizuo Kakutani.
\newblock On equivalence of infinite product measures.
\newblock {\em Ann. of Math. (2)}, 49:214--224, 1948.

\bibitem{MR2944069}
Shao-Gao Lv and Yun-Long Feng.
\newblock Integral {O}perator {A}pproach to {L}earning {T}heory with
  {U}nbounded {S}ampling.
\newblock {\em Complex Anal. Oper. Theory}, 6(3):533--548, 2012.

\bibitem{Maslowski}
{Maslowski, Bohdan and Nualart, David}.
\newblock {Evolution equation driven by a fractional Brownian motion}.
\newblock {\em {Journal of Functional Analysis}}, 202(1):277--305, 2003.

\bibitem{MR0214150}
Edward Nelson.
\newblock {\em Dynamical theories of {B}rownian motion}.
\newblock Princeton University Press, Princeton, N.J., 1967.

\bibitem{Nelson_flows}
Edward Nelson.
\newblock {\em Topics in dynamics. {I}: {F}lows}.
\newblock Mathematical Notes. Princeton University Press, Princeton, N.J.,
  1969.

\bibitem{MR2810909}
P.~Niyogi, S.~Smale, and S.~Weinberger.
\newblock A topological view of unsupervised learning from noisy data.
\newblock {\em SIAM J. Comput.}, 40(3):646--663, 2011.

\bibitem{MR0622034}
K.~R. Parthasarathy and K.~Schmidt.
\newblock {\em Positive definite kernels, continuous tensor products, and
  central limit theorems of probability theory}.
\newblock Springer-Verlag, Berlin, 1972.
\newblock Lecture Notes in Mathematics, Vol. 272.

\bibitem{MR1785620}
Yuval Peres, Wilhelm Schlag, and Boris Solomyak.
\newblock Sixty years of {B}ernoulli convolutions.
\newblock In {\em Fractal geometry and stochastics, {II}
  ({G}reifswald/{K}oserow, 1998)}, volume~46 of {\em Progr. Probab.}, pages
  39--65. Birkh\"auser, Basel, 2000.

\bibitem{MR1386842}
Yuval Peres and Boris Solomyak.
\newblock Absolute continuity of {B}ernoulli convolutions, a simple proof.
\newblock {\em Math. Res. Lett.}, 3(2):231--239, 1996.

\bibitem{MR751959}
Michael Reed and Barry Simon.
\newblock {\em Methods of modern mathematical physics. {I}}.
\newblock Academic Press Inc. [Harcourt Brace Jovanovich Publishers], New York,
  second edition, 1980.
\newblock Functional analysis.

\bibitem{saitoh}
S.~Saitoh.
\newblock {\em Theory of reproducing kernels and its applications}, volume 189.
\newblock Longman scientific and technical, 1988.

\bibitem{MR0102759}
I.~E. Segal.
\newblock Distributions in {H}ilbert space and canonical systems of operators.
\newblock {\em Trans. Amer. Math. Soc.}, 88:12--41, 1958.

\bibitem{MR0264761}
David Shale and W.~Forrest Stinespring.
\newblock Wiener processes. {II}.
\newblock {\em J. Functional Analysis}, 5:334--353, 1970.

\bibitem{MR2875207}
Guangjun Shen and Chao Chen.
\newblock Stochastic integration with respect to the sub-fractional {B}rownian
  motion with {$H\in(0,\frac12)$}.
\newblock {\em Statist. Probab. Lett.}, 82(2):240--251, 2012.

\bibitem{MR2488871}
Steve Smale and Ding-Xuan Zhou.
\newblock Online learning with {M}arkov sampling.
\newblock {\em Anal. Appl. (Singap.)}, 7(1):87--113, 2009.

\bibitem{MR1557013}
J.~von Neumann.
\newblock On infinite direct products.
\newblock {\em Compositio Math.}, 6:1--77, 1939.

\bibitem{MR1822853}
Yang Wang.
\newblock On the number of {D}aubechies scaling functions and a conjecture of
  {C}hyzak et al.
\newblock {\em Experiment. Math.}, 10(1):87--89, 2001.

\end{thebibliography}
\def\cprime{$'$} \def\lfhook#1{\setbox0=\hbox{#1}{\ooalign{\hidewidth
  \lower1.5ex\hbox{'}\hidewidth\crcr\unhbox0}}} \def\cprime{$'$}
  \def\cprime{$'$} \def\cprime{$'$} \def\cprime{$'$} \def\cprime{$'$}

\end{document}